\documentclass[10pt]{article}
\usepackage[active]{srcltx}

\usepackage[colorlinks=true]{hyperref}
\hypersetup{urlcolor=blue, citecolor=red}

\usepackage{amssymb,amsmath,amsfonts,bbm,pifont,upgreek,bbold,accents,ulem}

\font\teneufm=eufm10
\font\seveneufm=eufm7
\font\fiveeufm=eufm5
\newfam\eufmfam
\textfont\eufmfam=\teneufm
\scriptfont\eufmfam=\seveneufm
\scriptscriptfont\eufmfam=\fiveeufm
\def\eufm#1{{\fam\eufmfam\relax#1}}

\newcommand\beq[1]{ \begin{equation}\label{#1} }
\newcommand{\eeq}{ \end{equation} }
\newcommand{\beqno}{ \[ }
\newcommand{\eeqno}{ \] }
\newcommand\beqa[1]{ \begin{eqnarray} \label{#1}}
\newcommand{\eeqa}{ \end{eqnarray} }
\newcommand{\beqano}{ \begin{eqnarray*} }
\newcommand{\eeqano}{ \end{eqnarray*} }
\newcommand\arr[1]{\left\{\begin{array}{l}#1\end{array}\right.}
\renewcommand{\theequation}{\arabic{section}.\arabic{equation}}

\newtheorem{theorem}{Theorem}[section]
\newtheorem{definition}{Definition}[section]
\newtheorem{proposition}{Proposition}[section]
\newtheorem{lemma}{Lemma}[section]
\newtheorem{sublemma}{Sublemma}[section]
\newtheorem{remark}{Remark}[section]
\newtheorem{notationalremark}{Notational Remark}[section]
\newtheorem{corollary}{Corollary}[section]
\newtheorem{assumption}{Assumption}[section]
\newtheorem{claim}{Claim}[section]
\newtheorem{tools}{$\negsp\negsp$}[subsection]

\newcommand\thm[1]{ \begin{theorem}\label{#1}}
\newcommand\thmtwo[2]{ \begin{theorem}[#1]\label{#2}}
\newcommand\ethm{ \end{theorem} }
\newcommand\dfn[1]{ \begin{definition}\label{#1} \rm}
\newcommand\dfntwo[2]{ \begin{definition}[#1]\label{#2} \rm}
\newcommand\edfn{ \end{definition} }
\newcommand\pro[1]{ \begin{proposition}\label{#1}}
\newcommand\protwo[2]{ \begin{proposition}[#1]\label{#2}}
\newcommand\epro{ \end{proposition} }
\newcommand\lem[1]{ \begin{lemma}\label{#1}}
\newcommand\lemtwo[2]{ \begin{lemma}[#1]\label{#2}}
\newcommand\elem{ \end{lemma} }
\newcommand\sublem[1]{ \begin{sublemma}\label{#1}}
\newcommand\sublemtwo[2]{ \begin{sublemma}[#1]\label{#2}}
\newcommand\esublem{ \end{sublemma} }
\newcommand\rem[1]{ \begin{remark}\label{#1} \rm}
\newcommand\erem{ \end{remark} }
\newcommand\notrem[1]{ \begin{notationalremark}\label{#1} \rm}
\newcommand\enotrem{ \end{notationalremark} }
\newcommand\cor[1]{ \begin{corollary}\label{#1}}
\newcommand\cortwo[2]{ \begin{corollary}[#1]\label{#2}}
\newcommand\ecor{ \end{corollary} }
\newcommand\asmp[1]{ \begin{assumption}\label{#1}}
\newcommand\asmptwo[2]{ \begin{assumption}[#1]\label{#2}}
\newcommand\easmp{ \end{assumption} }
\newcommand\clm[1]{ \begin{claim}\label{#1}}
\newcommand\eclm{ \end{claim} }

\newcommand\equ[1]{{\rm (\ref{#1})}}

\expandafter\chardef\csname pre amssym.def
at\endcsname=\the\catcode`\@
\catcode`\@=11
\def\undefine#1{\let#1\undefined}
\def\newsymbol#1#2#3#4#5{\let\next@\relax
 \ifnum#2=\@ne\let\next@\msafam@\else
 \ifnum#2=\tw@\let\next@\msbfam@\fi\fi
 \mathchardef#1="#3\next@#4#5}
\def\mathhexbox@#1#2#3{\relax
 \ifmmode\mathpalette{}{\m@th\mathchar"#1#2#3}%
 \else\leavevmode\hbox{$\m@th\mathchar"#1#2#3$}\fi}
\def\hexnumber@#1{\ifcase#1 0\or 1\or 2\or 3\or 4\or 5\or 6\or 7\or
8\or
 9\or A\or B\or C\or D\or E\or F\fi}
\ifcase\@ptsize
 \font\tenmsb=msbm10
 \font\sevenmsb=msbm7
 \font\fivemsb=msbm5
\or
 \font\tenmsb=msbm10 scaled \magstephalf
 \font\sevenmsb=msbm7 scaled \magstephalf
 \font\fivemsb=msbm5  scaled \magstephalf
\or
 \font\tenmsb=msbm10 scaled \magstep1
 \font\sevenmsb=msbm7 scaled \magstep1
 \font\fivemsb=msbm5 scaled \magstep1
\fi
\newfam\msbfam
\textfont\msbfam=\tenmsb
\scriptfont\msbfam=\sevenmsb
\scriptscriptfont\msbfam=\fivemsb
\edef\msbfam@{\hexnumber@\msbfam}
\def\Bbb#1{\fam\msbfam\relax#1}
\def\widehat#1{\setboxz@h{$\m@th#1$}%
 \ifdim\wdz@>\tw@ em\mathaccent"0\msbfam@5B{#1}%
 \else\mathaccent"0362{#1}\fi}
\def\widetilde#1{\setboxz@h{$\m@th#1$}%
 \ifdim\wdz@>\tw@ em\mathaccent"0\msbfam@5D{#1}%
 \else\mathaccent"0365{#1}\fi}

\def\RIfM@{\relax\ifmmode}
\def\nonmatherr@#1{\errmessage{\string#1\space allowed only in math mode}}
\def\Bbb{\RIfM@\expandafter\Bbb@\else
 \expandafter\nonmatherr@\expandafter\Bbb\fi}
\def\Bbb@#1{{\Bbb@@{#1}}}
\def\Bbb@@#1{\fam\msbfam\relax#1}
\def\setboxz@h{\setbox\z@\hbox}
\def\wdz@{\wd\z@}
\catcode`\@=\csname pre amssym.def at\endcsname
\newcommand{\ie}{{\rm i.e.\,}}
\newcommand{\eg}{{\it e.g.\,}}

\newcommand{\nl}{{\smallskip\noindent}}
\newcommand{\noi}{{\noindent}}

\newcommand{\negsp}{\hspace{-.09truecm}}  

\newcommand{\dst}{\displaystyle}

\newcommand{\torus}{ {\Bbb T}   }

\newcommand{\real}{ {\Bbb R}   }
\newcommand{\integer}{ {\Bbb Z}   }

\renewcommand{\a }{ {\alpha}   }
\renewcommand{\b}{ {\beta}   }
\newcommand{\g}{ {\gamma}   }
\newcommand{\G}{ {\Gamma}   }
\renewcommand{\d}{ {\delta}   }

\renewcommand{\k}{ {\kappa}   }
\renewcommand{\l}{ {\lambda}   }
\renewcommand{\L}{ {\Lambda}   }
\newcommand{\m}{ {\mu}   }
\newcommand{\n}{ {\nu}   }

\newcommand{\p}{ {\pi}   }

\renewcommand{\t}{ {\tau}   }

\renewcommand{\O}{ {\Omega}   }

\newcommand{\const}{{\, \rm const\, }}

\newcommand{\cA}{ {\cal A} }

\newcommand{\cR}{ {\cal R} }

\newcommand{\cM}{ {\cal M} }

\newcommand\ppu{{ (1) }}
\newcommand\ppd{{ (2) }}
\newcommand\ppt{{ (3) }}

\newcommand\ppj{{ (j) }}

\newcommand\ppn{{ (n) }}

\newcommand\ppi{{ (i) }}

\newcommand\ppo{{ (0) }}

\newcommand\ul{{\uplambda}}
\newcommand\ux{{\upxi}}
\newcommand\uh{{\upeta}}
\newcommand\up{{\rm p}}
\newcommand\uq{{\rm q}}
\newcommand\uz{{\rm z}}

\newcommand\id{{\, \rm id \,}}

\begin{document}

\title{Canonical coordinates for the planetary problem\footnote{ Research supported by ERC Ideas-Project 306414 ``Hamiltonian PDEs and small divisor problems: a dynamical systems approach'' and STAR Project of Federico II University, Naples.}}

\author{  
Gabriella Pinzari\\
\vspace{-.2truecm}
{\footnotesize Dipartimento di Matematica ed Applicazioni ``R. Caccioppoli''}
\\{\footnotesize Universit\`a di Napoli ``Federico II"}
\vspace{-.2truecm}
\\{\footnotesize Monte Sant'Angelo -- Via Cinthia I-80126 Napoli (Italy)}
\vspace{-.2truecm}
\\{\scriptsize gabriella.pinzari@unina.it}
}\date{November, 3, 2014}
\maketitle
\vskip.1in
\noi

\begin{abstract}\footnotesize{We collect classical  and more recent material  about possible  symplectic descriptions  of the phase space of the planetary problem. 
}
\end{abstract}

{
\nl
{\footnotesize{\bf Keywords:} Canonical coordinates.  Planetary problem. Delaunay--Poincar\'e coordinates. Jacobi's reduction of the nodes. Deprit's reduction of the nodes. RPS variables.  Keplerian ellipses. Symmetries.  Perihelia reduction. 

\nl
}

\maketitle

\tableofcontents

\section{Introduction}
\setcounter{equation}{0}
\label{intro}
\renewcommand{\theequation}{\arabic{equation}}
\paragraph{\bf 1.1} In 1963, V. I. Arnold stated his celebrated Theorem\footnote{{\cite[Chapter III, p. 125]{arnold63}}{
\sl ``For the majority of initial conditions under which the instantaneous orbits of the planets are close to circles lying in a single plane, perturbation of the planets on one another produces, in the course of an infinite interval of time, little change on these orbits provided the masses of the planets are sufficiently small. {\rm[...]}
In particular {\rm[...]} in the n-body problem there exists a set of initial conditions having a positive Lebesgue measure and such that, if the initial positions and velocities of the bodies belong to this set, the distances of the bodies from each other will remain
perpetually bounded.''}}  on  the Stability of Planetary Motions (hereafter, Arnold's Theorem). The general  proof of his wonderful statement (that he provided completely only for the particular case of three bodies constrained on a plane) turned out to be more difficult than expected and was next completed by  J. Laskar, P. Robutel, M. Herman, J. F\'ejoz, L. Chierchia and the author. We refer the reader to the technical papers \cite{arnold63},  \cite{laskarR95},  \cite{robutel95},  \cite{maligeRL02},  \cite{herman09},  \cite{fejoz04},  \cite{pinzari-th09},  \cite{chierchiaPi11b} for detailed information; to \cite{fejoz13},  \cite{chierchia13},  \cite{chierchiaPi14}, or  the introduction of \cite{pinzari13} for  reviews.

\nl
The complete understanding of Arnold's Theorem relied on an analytic part and a geometric one, both  highly non trivial.  Of such two aspects,  the analytic part  was basically  settled out  since \cite{arnold63} (notwithstanding refinements next given in \cite{fejoz04},  \cite{chierchiaPi10}). The geometrical aspects, of which Arnold just provided, for the general case, only sketchy conjectures, were instead mostly  unexplored after his 1963's paper and  have been  only recently  clarified \cite{pinzari-th09},  \cite{chierchiaPi11b}, \cite{pinzari14}.

\nl
In fact, switching from the three--body case to the many--body one needed to develop new constructions  not known at those times.

\nl
The purpose of this note is to provide a historical survey 
of this  latter part. We shall describe previous classical approaches
going back to Delaunay, Poincar\'e, Jacobi and point out more recent progresses, based on the papers \cite{pinzari-th09},  \cite{chierchiaPi11a},  \cite{chierchiaPi11b},  \cite{chierchiaPi11c},  \cite{pinzari13},  \cite{pinzari14}.

\nl
We begin with introducing the problem.

\vskip.1in
\noi
The planetary problem consists  in studying the motion of $(1+n)$ point masses, a ``sun'' and $n$ ``planets''  interacting through gravity. 
This is a Hamiltonian problem: If the masses of the sun and of the planets are denoted, respectively, as $m_0$, $\m m_1$, $\cdots$, $\m m_n$, where  $\m$ is a very small number, and the Euclidean length as $|\, \cdot\, |$, the motion equations 
may be written in the form of Hamilton equation
where the Hamiltonian (the energy\footnote{More precisely, $\m{\rm H}_{(1+n)\rm b}(p/\m,q)
$
 corresponds to be the mechanical energy $T+U$.
 }) of the system is
 $${\rm H}_{(1+n)\rm b}(p,q)=\sum_{1\le i\le n}\frac{|p^\ppi|^2}{2 m_i}-
\sum_{1\le i\le n}\frac{m_0m_i}{|q^\ppo-q^\ppi|}+\m\frac{|p^\ppo|^2}{2m_0}-\m\sum_{1\le i<j\le n}\frac{m_im_j}{|q^\ppi-q^\ppj|}\ .$$
Here, for $0\le i\le n$, $1\le j\le 3$, $q=(q^\ppo, \cdots, q^\ppn)$, with  $q^\ppi=(q^\ppi_1, q^\ppi_2, q^\ppi_3)$ are the planets' positions   and $p=(p^\ppo, \cdots, p^\ppn)$,  with $p^\ppi=(p^\ppi_1, p^\ppi_2, p^\ppi_3)$, their conjugated momenta.
 
 \nl
Despite of its apparent  $(3+3n)$ degrees of freedom 
 ($(2n+2)$ for the problem in the plane), ${\rm H}_{(1+n)\rm b}$ possesses  a variety of  integrals
. These are: the three components of the center of mass 
 ${\rm Q}=(m_0q^\ppo+\m\sum_{1\le i\le n}m_i q^\ppi)/\rm const$; the three components of the total linear momentum
 ${\rm P}=\sum_{0\le i\le n}p^\ppi=\const\dot{\rm Q}$; the three components of the total angular momentum
 ${\rm C}=\sum_{0\le i\le n}q^\ppi\times p^\ppi$ and the energy ${\rm E}={\rm H}_{(1+n)\rm b}$.
 %
 %
 %
 %
Here, ``$\const$'' is the sum of the masses and ``$\times$'' denotes skew--product.

\nl
Poincar\'e proved that no other analytic integral for the system may be found, apart for the ten ones listed above. They are more than enough to integrate the problem ``by quadratures'' in the case of two bodies ($n=1$; see \S \ref{2 bodies}). When $n\ge 2$, the problem is non--integrable. Its dynamics may be very complicated and exhibit also chaotic behaviors \cite{fejozGKR14}. 

\nl
The integrals above are  independent but do not commute. It 
 is however possible to isolate among them  five independent  and pairwise commuting quantities, for the problem in the space; three of them for the problem constrained on a plane. These are:   the three components of the linear momentum ${\rm P}$, the third component ${\rm C}_3$ of the total angular momentum ${\rm C}$ and\footnote{Recall that the three components of ${\rm C}$ do not commute. Rather, they verify $\{{\rm C}_i, {\rm C}_j\}=\d_{ijk}{\rm C}_k$, where $\d_{ijk}$ is the Levi--Civita symbol.} its Euclidean length ${\rm G}:=|{\rm C}|$ in the former case; the two components of ${\rm P}$ and ${\rm C}_3$, in the latter. This  tells us that the ``effective'' number of degrees  of freedom is   $(3n+3)-5=3n-2$ in the space; $(2n+2)-3=2n-1$ in the plane.

 \nl
 To any of such independent  integrals one  associates  (by Noether theorem) a one--parameter group of  symplectic transformations which leave ${\rm H}_{(1+n)\rm b}$  unvaried. For problem in the space, these are

 \begin{itemize}
 \item[\tiny\textbullet] translations $$\dst {\rm t}_a:\ (p^\ppi, q^\ppi)\to (p^\ppi, q^\ppi+a)\qquad a\in \real^3$$  with respect to three independent directions  $a=(a_1,a_2,a_3)$;
 
 \item[\tiny\textbullet]   ``synchronous'' rotations for the $p$ and the $q$--variables  about the $k^\ppt$--axis:
 $$\dst {\rm R}_{\a}^\ppt:\  (p^\ppi, q^\ppi)\to ({\rm R}^\ppt_\a p^\ppi, {\rm R}^\ppt_\a q^\ppi)\qquad \a\in \torus\ ;$$

 \item[\tiny\textbullet] or  about the ${\rm C}$--axis:
 $$\dst {\rm R}^{{\rm C}}_\b:\  (p^\ppi, q^\ppi)\to ({\rm R}^{\rm C}_{\b}p^\ppi, {\rm R}^{\rm C}_{\b}q^\ppi)\qquad \b\in \torus\ .$$

\nl
 Here,  ${\rm F}_0=(k^\ppu, k^\ppd, k^\ppt)$ is a prefixed reference frame, 
 $\torus:=\real/2\p\integer$ is the flat torus and ${\rm R}^\ppt_\a$, ${\rm R}^{\rm C}_\b$ are suitable $3\times 3$ matrices having the form of ${\rm R}^\ppt({\rm h})$ below (see Eq. \equ{R1R3}) in suitable reference frames, with ${\rm h}=\a, \b$.
 \end{itemize}
 
 \nl
  For the problem in the plane one has analogous symmetries. 
    

\nl 
We remark that, since the transformations listed above are caused by integrals, unless such integrals are eliminated in some way,  they are present in the Hamiltonian (in different forms), no matter what system of coordinates is used.
Arnold in 1963 argued that the symmetry by  rotations  
 plays a key r\^ole in order to overcome the problem of the ``proper degeneracy'' (that we shall recall in \S \ref{2 bodies}) of the planetary system, since it implies the existence of an elliptic equilibrium point for the secular\footnote{In this context, ``secular system'' is used to denote a suitable time--average of it. A rigorous definition will be given in \S \ref{Kepler maps*} below.}  system written in Poincar\'e coordinates (see \S \ref{2 bodies}) in correspondence of circular, planar motions. 
He foresaw this equilibrium point might be  a bifurcation  point  of full--dimensional invariant tori surrounded of quasi--periodic motions. Quite paradoxically,  the invariance by rotations itself  is just the origin of a strong degeneracy, for the problem in the space, that, if not suitably treated, strongly prevents the direct application of his project.
  Arnold's program was then successfully completed in \cite{robutel95} and, in its full generality, in \cite{pinzari-th09}, \cite{chierchiaPi11b}.

\nl
It is worth  to recall also that further symmetries for ${\rm H}_{(1+n)\rm b}$ do exist for the planetary Hamiltonian which are not related to integrals.

\nl
These  are 
  \begin{itemize}
  \item[\tiny\textbullet] ``asynchronous'' rotations for the $p$ and the $q$--variables, \ie,
   $$\dst (p^\ppi, q^\ppi)\to ({\rm R} p^\ppi, {\rm S} q^\ppi)$$
where ${\rm R}{\rm R}^{\rm t}=\id={\rm S}{\rm S}^{\rm t}$, with ${}^{\rm t}$ denoting transpose and ${\rm R}\ne{\rm S}$;

  \item[\tiny\textbullet] ``reflections'', \ie,  transformations of the form
\beq{reflections} \big((y^\ppi_1, y^\ppi_2, y^\ppi_3), (x^\ppi_1, x^\ppi_2, x^\ppi_3)\big)\to \big((r_1y^\ppi_1, r_2y^\ppi_2, r_3y^\ppi_3), (s_1x^\ppi_1, s_2x^\ppi_2, s_3x^\ppi_3)\big)\eeq
  with $r_j$, $s_j=\pm1$.
 \end{itemize}

\nl
These symmetries play an important r\^ole  in the case  one would like to know  whether it is possible to reduce completely the number of degrees of freedom of the system and, simultaneously, keep some parity (compare also the next paragraph).

\paragraph{\bf 1.2} In general for a quasi--integrable system
$${\rm H}(I,\varphi)=h(I)+\varepsilon f(I,\varphi)\qquad (I,\varphi)\in V\times \torus^n\qquad V\subset \real^n\quad{\rm open}$$
``natural'' action--angle  coordinates $(I,\varphi)$
 are  uniquely determined (up to isomorphisms) by the integrable part $h$.

\nl
In the case of the planetary problem, notwithstanding its clear quasi--integrable structure (it is close to the integrable problem consisting of  the uncoupled interaction of each planet with the sun), its proper degeneracies mentioned above   imply that  the choice of symplectic coordinates  is not uniquely determined. 
We then aim to describe   possible 
sets of canonical coordinates which let the system free of  its integrals (all of, or just a part of them). We shall only deal with systems of coordinates such that the unperturbed part takes the aspect of the classical Keplerian form (see Eq. \equ{Kepler} below). See also \cite{whittaker59} and references therein for further examples.



\nl
In order to reduced symmetries due to integrals, one has to deal with the linear momentum ${\rm P}$ and   the angular momentum ${\rm C}$.\\
The reduction of the linear momentum may be performed in two (equivalent)  ways, as described in \S \ref{Linear Momentum Reductions}, switching to {\sl Jacobi} or {\sl heliocentric}  coordinates. \\
The  reduction of the angular momentum does not present difficulties in the case of the planar problem and can be performed using common tools of canonical transformations. Therefore, we shall not deal about it (see however \cite{chierchiaPi11c} for recent results concerning the Birkhoff normal form for the completely reduced problem in the planar case). \\
The  reduction of the angular momentum for the problem in the space is instead something more intriguing.\\
For a long time, only the spatial three--body case could be handled, by means of a tool, known since the XIX century, under the name of   ``reduction of the nodes'', developed   by Jacobi and Radau \cite{jacobi1842}, \cite{radau1868}. It has been  used in \cite{robutel95} to prove Arnold's Theorem in this case and is recalled in \S \ref{Jacobi Radau}.\\
In lack of coordinates fitted to reduce completely rotations for the spatial problem in the case of more than three bodies, many authors \cite{abdullahA01}, \cite{herman09}, \cite{fejoz04}, \cite{nehorosev77} used 
classical  sets coordinates named after Delaunay and Poincar\'e. Delaunay and Poincar\'e coordinates (which are recalled in \S \ref{2 bodies}) are suited to the quasi--integrability of the problem, but  not  to its symmetry by rotations. \\
Only recently, new systems of coordinates have been applied to the problem \cite{pinzari-th09}, \cite{chierchiaPi11b}, in order to prove {\sl directly} Arnold's Theorem and estimating the measure of its stable motions, in terms of the maximum of eccentricities and inclinations of unperturbed Keplerian motions. These coordinates are related to a certain set discovered in the 80's by  A. Deprit \cite{deprit83}, and re-discovered by the author during her PhD \cite{pinzari-th09}, which reduce {\sl completely} the number of degrees of freedom of the system. 
Deprit's coordinates generalize to an arbitrary number of planets the  reduction of the nodes by Jacobi and Radau, and, as well as it,  are defined only for the problem in the space. 
To overcome the problem of singularity of the coordinates for planar motions, in \cite{pinzari-th09}, \cite{chierchiaPi11b} (see also a conjecture in \cite{arnold63} and a formal construction in \cite{maligeRL02}),
starting from Deprit's coordinates,  a system of coordinates realizing just a {\sl partial} reduction has been introduced. This allows to let the system free from its rotational degeneracy mentioned above (and hence to obtain the complete proof of Arnold's theorem), even though with an extra degree of freedom.  With this extra--degree of freedom it is also possible to treat the spatial and the planar problem in a unified way:  the partially reduced coordinates of the spatial case are defined also in the planar case, where they reduce to the Poincar\'e (hence, unreduced) planar coordinates. This part of the story is told  in the following \S \ref{Deprit' s reduction}--\S \ref{regularizations}.\\
We now remark that, as it follows from  the papers \cite{chierchiaPi11b} and \cite{chierchiaPi11c}, {\sl complete}, global reductions of the number of degrees of freedom are available  both for  the planar and the spatial problem {\sl separately}. This has the following physical consequence:  almost co-planar  quasi--periodic motions of the spatial problem with the minimum number of degrees of freedom need not to be closer and closer to the corresponding motions of the planar problem. Clearly, one should reasonably expect the opposite situation. It is to add that a certain continuity between the spatial and the planar three body problem was invoked by Arnold since \cite{arnold63} and this turned out to be a controversial question widely discussed since M. Herman's  investigations to this problem. In  \S \ref{Full reduction and  symmetries}--\S\ref{perihelia reduction} we  discuss the problem of having a unique, well defined, global  system of coordinates for the spatial and the planar problem. We shall present a  system of coordinates (which we name ``perihelia reduction'') which does this job and, moreover,  keeps parities due to reflections. This new set  reveals to be useful  to prove (besides the aforementioned continuity) that 
 quasi--periodic motions with minimum number of independent frequencies do exist even away from Arnold's elliptic equilibrium, \ie, for relatively large eccentricities and inclinations. See \cite{pinzari14} for an announcement of this result.

\section{Linear momentum reductions}\label{Linear Momentum Reductions}

The first step  consists into eliminating the motion of the sun by fixing  the coordinates  of the center of mass ${\rm Q}$. This leads to a $3n$--degrees of freedom Hamiltonian which governs the motion of the planets. The methods used in literature to achieve  this reduction are essentially two, respectively referred to as the {\sl heliocentric reduction} and the reduction via {\sl Jacobi coordinates}.

\nl
The heliocentric reduction has been used in \cite{herman09}, \cite{fejoz04}, \cite{robutel95}, \cite{chierchiaPi11b} and goes as follows. 
One performs the change of variables
 $$
\arr{\dst x^\ppo:=q^\ppo\\ \\
x^\ppi:=q^\ppi-q^\ppo\quad 1\le i\le n
}\qquad\qquad \arr{\dst y^\ppo:=\sum_{i=0}^n p^\ppi={\rm P}\\ \\
y^\ppi:=p^\ppi\quad 1\le i\le n\ .
}
 $$
 The change is linear, canonical (more precisely, it is homogeneous\footnote{As usual, by ``canonical, we mean that the standard 2--form $\sum_{0\le i\le n}\sum_{1\le j\le 3}dp^\ppi_j\wedge dq^\ppi_j$ is preserved. By ``homogeneous'', we mean that the standard 1--form $\sum_{0\le i\le n}\sum_{1\le j\le 3}p^\ppi_jdq^\ppi_j$ is preserved.}). The conservation of $ y^\ppo={\rm P}$ implies $x^\ppo$ is cyclic and, on the manifold where ${\rm Q}$ is constant, we have ${\rm P}=0$. Therefore, we can conveniently fix  $ x^\ppo=0$ and ${\rm P}=0$ in order to find the expression of the reduced Hamiltonian. This amounts to take
 $$
\arr{\dst q^\ppo=0\\ \\
q^\ppi=x^\ppi\quad 1\le i\le n
}\qquad\qquad \arr{\dst p^\ppo=-\sum_{i=1}^n y^\ppi\\ \\
p^\ppi=y^\ppi\quad 1\le i\le n\ .
}
 $$
Correspondingly, we find for the ``heliocentric Hamiltonian'' the expression
\beqa{helio}
{\rm H}_{\rm hel}(y, x)
&=&\sum_{i=1}^n
\big(\frac{|y^\ppi|^2}{2{\eufm m}_i}-\frac{{\eufm m}_i{\eufm M}_i}{|x^\ppi|}\big)+\m\sum_{1\le i<j\le n}\big(\frac{y^\ppi\cdot y^\ppj}{m_0}-\frac{m_im_j}{|x^\ppi-x^\ppj|}\big)\nonumber\\
&=:&\sum_{i=1}^n
\big(\frac{|y^\ppi|^2}{2{\eufm m}_i}-\frac{{\eufm m}_i{\eufm M}_i}{|x^\ppi|}\big)+\m f_{\rm hel}(y, x)
\eeqa
where $\dst {\eufm m}_i:=\frac{m_0 m_i}{m_0+\m m_i}$ and ${\eufm M}_i:={m_0+\m m_i}$.

\vskip.1in
\noi
The reduction via Jacobi coordinates has been used, for example, in \cite{arnold63}, \cite{zhao14}. It works well  when the number of planets is small. For simplicity, we present it for the three--body case; generalizations may be obtained.

\nl
Consider three masses, $m_0$, $\m m_1$ and $\m m_2$. Let
$${\rm Q}':=\frac{m_0q^\ppo+\m m_1 q^\ppu}{m_0+\m m_1}
$$
 the center of mass of $m_0$ and $\m m_1$ 
and switch to the (homogeneous--canonical) coordinates
 $$\arr{
\dst  \tilde x^\ppo=q^\ppo\\ \\
\dst \tilde x^\ppu=q^\ppu-q^\ppo\\ \\
\dst \tilde x^\ppd=q^\ppd-{\rm Q}'
 }\qquad \arr{
\dst  \tilde y^\ppo=p^\ppo+p^\ppu+p^\ppd\\ \\
\dst \tilde y^\ppu=p^\ppu+\frac{\m m_1}{m_0+\m m_1}p^\ppd\\ \\
\dst \tilde y^\ppd=p^\ppd\ .
 }$$
Proceeding in a similar way as for the heliocentric reduction (\ie, using the cyclic character of $\widetilde x^\ppo$ and restricting to the manifold where ${\rm Q}$ is constant, which allows to take also $\tilde y^\ppo=0$) the expression of the transformed Hamiltonian can be found taking
\beq{Jac}\arr{
\dst  q^\ppo=0\\ \\
\dst q^\ppu=\tilde x^\ppu\\ \\
\dst q^\ppd=\frac{\m m_1}{m_0+\m m_1}\tilde x^\ppu+\tilde x^\ppd
 }\qquad \arr{
\dst  p^\ppo=-\tilde y^\ppu-\frac{m_0}{m_0+\m m_1}\tilde y^\ppd\\ \\
\dst p^\ppu=\tilde y^\ppu-\frac{\m m_1}{m_0+\m m_1}\tilde y^\ppd\\ \\
\dst p^\ppd=\tilde y^\ppd\ .
 }\eeq
Correspondingly, one finds the ``Jacobi--Hamiltonian''

\beqa{Jac*}{\rm H}_{\rm Jac}&=&
\frac{|\tilde y^\ppu|^2}{2\tilde {\eufm m}_1}+\frac{|\tilde y^\ppd|^2}{2\tilde{\eufm m}_2}
-
\frac{m_0m_1}{|\tilde x^\ppu|}-
\frac{m_0m_2}{|\dst\tilde x^\ppd+\frac{\m m_1}{m_0+\m m_1}\tilde x^\ppu|}\nonumber\\
&-&\m\frac{m_1m_2}{\dst|\frac{m_0}{m_0+\m m_1}\tilde x^\ppu-\tilde x^\ppd|}
\eeqa
that we write in the form
\beqano{\rm H}_{\rm Jac}
&=:&\sum_{i=1}^2
\big(\frac{|\tilde y^\ppi|^2}{2\tilde{\eufm m}_i}-\frac{\tilde{\eufm m}_i\tilde{\eufm M}_i}{|x^\ppi|}\big)+\m f_{\rm Jac}(\tilde y, \tilde x)
\eeqano
with $\dst 
\tilde {\eufm m}_1={\eufm m}_1=\frac{m_0m_1}{m_0+\m m_1}$, $\dst\tilde{\eufm m}_2=\frac{m_2(m_0+\m m_1)}{m_0+\m m_1+\m m_2}={\eufm m}_2+{\rm O}(\m^2)$, $\tilde{\eufm M}_1:={\eufm M}_1=m_0+\m m_1$, $\dst\tilde{\eufm M}_2:=m_0\frac{m_0+\m m_1+\m m_2}{m_0+\m m_1}$.

\nl
The choice of $\tilde{\eufm M}_i$  is justified by the fact that, if one takes the expansion of the two last terms of ${\rm H}_{\rm Jac}$ in powers of $\m$, then he finds the Hamiltonian
\beq{Jac**}\widetilde{\rm H}_{\rm Jac}=
\sum_{1\le i\le 2}\big(\frac{|\tilde y^\ppi|^2}{2\tilde {\eufm m}_i}-
\frac{\tilde{\eufm m}_i\tilde{\eufm M}_i}{|\tilde x^\ppi|}\big)+\m m_1m_2\big(\frac{\tilde x^\ppu\cdot \tilde x^\ppd}{|\tilde x^\ppd|^3}
-\frac{1}{\dst|\tilde x^\ppu-\tilde x^\ppd|}\big)\ ,\eeq
which differs from ${\rm H}_{\rm Jac}$ by ${\rm O}(\m^2)$. 
The Hamiltonian $\widetilde{\rm H}_{\rm Jac}$ has been proposed in \cite{arnold63} for the proof of Arnold's theorem.

\vskip.1in
\nl
The relation between the heliocentric and the Jacobi reduction is simple.

\nl
The definitions above imply that, if $(y^\ppu$, $y^\ppd$, $x^\ppu$, $x^\ppd)$ $=$ $\phi_{\rm hel/Jac}$ $(\tilde y^\ppu$, $\tilde y^\ppd$, $\tilde x^\ppu$, $\tilde x^\ppd)$ denotes the linear, symplectic, close to the identity transformation given by 

\beqno\arr{
\dst x^\ppu=\tilde x^\ppu\\ \\
\dst x^\ppd=\frac{\m m_1}{m_0+\m m_1}\tilde x^\ppu+\tilde x^\ppd
 }\qquad \arr{
\dst y^\ppu=\tilde y^\ppu-\frac{\m m_1}{m_0+\m m_1}\tilde y^\ppd\\ \\
\dst y^\ppd=\tilde y^\ppd
 }\eeqno
 then the heliocentric and the Jacobi Hamiltonian are related by
\beq{hel/Jac}{\rm H}_{\rm Jac}={\rm H}_{\rm hel}\circ\phi_{\rm hel/Jac}\ . \eeq
From this relation,  a certain ``equivalence'' (in the sense of normal form) between ${\rm H}_{\rm Jac}$, ${\rm H}_{\rm hel}$ follows, that will be discussed later (compare \S \ref{Kepler maps*}).

\section{Canonical coordinates fitted to the perturbative setting}
Whatever is the method which one chooses to eliminate the coordinates of the sun (heliocentric, Jacobi or others), the resulting  Hamiltonian  for the planets takes the form of a quasi--integrable system, where 
\begin{itemize}
\item[\tiny\textbullet]  the leading part (of ``order one'')  is a sum of uncoupled ``two--body systems'' of the form
 \beq{two bodies}\dst h^\ppi_{\rm2b}=\frac{|y^\ppi|^2}{2{\eufm m}_i}-\frac{{\eufm m}_i{\eufm M}_i}{|x^\ppi|}\ ,\eeq
which, as well known, are integrable; 

\item[\tiny\textbullet] the perturbing function (of ``order $\m$''), as well as $\dst h^\ppi_{\rm2b}$, possesses the integral
\beq{C}{\rm C}=\sum_{1\le i\le n}{\rm C}^\ppi\eeq
where ${\rm C}^\ppi=x^\ppi\times y^\ppi$.
\end{itemize}

\subsection{On the canonical integration of the two--body problem (Delaunay and Poincar\'e coordinates)}\label{2 bodies}

\nl
Referring to classical textbooks for more details, in this section we  recall a few facts about the construction of canonical variables for $h^\ppi_{\rm2b}$.

\nl
One often says that  $h^\ppi_{\rm 2b}$ is ``super--integrable''. This locution reflects the fact that $h^\ppi_{\rm 2b}$, despite of having three degrees of freedom,  possesses five independent integrals, even more  than its degrees of freedom.  These are:  the three components of the angular momentum ${\rm C}^\ppi=x^\ppi\times y^\ppi$, the direction of the so--called eccentricity vector ${\rm L}^\ppi$, which is perpendicular to ${\rm C}^\ppi$ and the energy. 
The presence of so many integrals causes two degeneracies  in  the integration of $h^\ppi_{\rm 2b}$, which are often referred to as ``proper degeneracies''. More precisely, degeneracies arise from the fact that, excluding the energy, the four remaining  integrals may be rearranged into two couples of canonical conjugated variables ($({\rm H}_i, {\rm h}_i)$ and $(\G_i, {\rm g}_i)$ below).  Let us recall this construction.

\nl
The integrals above are not in all  involution. Then one  chooses the Euclidean length $|{\rm C}^\ppi|$ and the third component ${\rm C}^\ppi_3={\rm C}^\ppi\cdot k^\ppt$ (where ${\rm F}_0=(k^\ppu, k^\ppd, k^\ppt)$ is a prefixed orthonormal triple in $\real^3$), which are in mutual involution and also commute with the energy, and then introduces a change of variables which has such functions among the generalized momenta, defined as follows.

 \nl
  Assume that ${\rm C}^\ppi$ is not parallel to $k^\ppt$ and let 
  \beq{D-nodes}n_i:=k^\ppt\times {\rm C}^\ppi\ ,\eeq 
so that $n_i\ne 0$  well defines the intersection (``node'') of the $(k^\ppu, k^\ppd)$ plane with the plane orthogonal to ${\rm C}^\ppi$ (the plane of the orbit\footnote{By the conservation of ${\rm C}^\ppi$, the orbit of $h^\ppi_{\rm 2b}$ lies in a plane. }).    Then define the following coordinates
\beq{reduct}
\arr{\dst {\rm H}_i={\rm C}^\ppi_3\\ \\
\dst {\rm R}_i:=\frac{y^\ppi\cdot x^\ppi}{|x^\ppi|}\\ \\
\dst\Phi_i:=|{\rm C}^\ppi|
}\qquad \arr{\dst {\rm h}_i=\a_{k^\ppt}(k^\ppu, n_i)\\ \\
\dst {\rm r}_i:=|x^\ppi|\\ \\
\dst \varphi_i:=\a_{{\rm C}^\ppi}(n_i, x^\ppi)\ .
}
\eeq
Here, given three independent vectors $u$, $v$, $w\in \real^3$, with $w\perp(u, v)$,  $\a_w(u, v)\in \torus$ denotes the oriented angle formed by $(u, v)$ in the positive (counterclockwise) verse determined by $w$.

\nl
The inverse formulae of \equ{reduct} are as follows. If
\beq{bar y bar x}\bar x^\ppi({\rm R}_i, \Phi_i, {\rm r}_i, \varphi_i):=\arr{\dst {\rm r}_i\cos\varphi_i\\ \\
\dst {\rm r}_i\sin\varphi_i
} \quad \bar y^\ppi({\rm R}_i, \Phi_i, {\rm r}_i, \varphi_i):=\arr{\dst {\rm R}_i\cos\varphi_i-\frac{\Phi_i}{{\rm r}_i}\sin\varphi_i\\
\dst {\rm R}_i\sin\varphi_i+\frac{\Phi_i}{{\rm r}_i}\cos\varphi_i
} \eeq
then
\beqa{x***}x^\ppi({\rm R}_i, \Phi_i, {\rm r}_i, \varphi_i)&=&{\rm R}^\ppt({\rm h}_i){\rm R}^\ppu(\iota_i) \bar x^\ppi({\rm R}_i, \Phi_i, {\rm r}_i, \varphi_i)\nonumber\\
 y^\ppi({\rm R}_i, \Phi_i, {\rm r}_i, \varphi_i)&=&{\rm R}^\ppt({\rm h}_i){\rm R}^\ppu(\iota_i) \bar y^\ppi({\rm R}_i, \Phi_i, {\rm r}_i, \varphi_i)
\eeqa
where $\iota_i$ denote the convex angles between 
${\rm C}$ and ${\rm C}^\ppi$, determined by
$$\cos\iota_i=\frac{{\rm H}_i}{\Phi_i}\ ;$$
 ${\rm R}^\ppu(\iota)$ ,
${\rm R}^\ppt({\rm h})$ are the elementary rotations
\beqa{R1R3} {\rm R}^\ppu(\iota)= \left(
 \begin{array}{ccc}
 1&0&0\\
0 &\cos i&-\sin i\\
0 &\sin i&\cos i \end{array}
 \right)
 \quad \quad {\rm R}^\ppt({\rm h})=\left(
 \begin{array}{ccc}
 \cos{\rm h}&-\sin{\rm h}&0\\
 \sin{\rm h}&\cos{\rm h}&0\\
 0&0&1
 \end{array}
 \right)\ .
 \eeqa

\nl
The change \equ{reduct} is homogeneous
and $h^\ppi_{\rm 2b}$ is transformed into
\beq{2bpc}h^\ppi_{\rm pc}:=\frac{{\rm R}_i^2}{2{\eufm m}_i}+\frac{\Phi_i^2}{2m_i{\rm r}_i^2}-\frac{{\eufm M}_i{\eufm m}_i}{{\rm r}_i}\ .\eeq
As expected, the angles  ${\rm h}_i$ and $\varphi_i$ disappear from $h^\ppi_{\rm pc}$ since they are conjugated  to the  integrals ${\rm H}_i$, $\Phi_i$.  The first degeneracy of the problem consists in the fact that the {\sl action} ${\rm H}_i$ disappears. 
 This happens because its conjugated angle ${\rm h}_i$ is {\sl also} an integral  (it is again related  to ${\rm C}^\ppi$) . 
 
\nl
We shall see in \S \ref{Deprit' s reduction}  that this ``rotational'' degeneracy of the two--body problem has a full generalization in the many--body case. This remarkable fact has been only recently pointed out \cite{pinzari-th09}, \cite{chierchiaPi11b}, \cite{chierchiaPi11c}.

\nl
The second step consists into the integration of the one--dimensional Hamiltonian $h^\ppi_{\rm pc}$, via the Hamilton--Jacobi method.

\nl
We are only interested to do it in the region of phase space where $h^\ppi_{\rm pc}<0$, where the energy levels $h^\ppi_{\rm 2b}={\rm E}^\ppi_{\rm 2b}$ are compact. In this case, one finds a two--dimensional canonical change variables parametrized by $\Phi_i$
\beq{lin Del}(\L_i,  \ell_i)\to ({\rm R}_i, {\rm r}_i)\eeq
which integrates $h^\ppi_{\rm pc}$. The integrated form of  $h^\ppi_{\rm pc}$ is just the ``Keplerian'' one
\beq{Kepler}h_{\rm K}^\ppi=-\frac{{\eufm m}_i^3{\eufm M}_i^2}{2\L_i^2}\ .\eeq
As expected, the angles $(\ell_i$, $ {\rm g}_i)$ disappear since they are conjugated to the integrals $(\L_i,\G_i)$. The fact that also $\G_i$ disappears (this is the second degeneracy of the problem) is related to the integral ${\rm L}_i$: the angle ${\rm g}_i$, conjugated to $\G_i$,  just corresponds to be the angle $\widehat{ n_i{\rm L}_i}$, and this angle does not move, so also $\G_i$ disappears. 

\nl
The motion then reduces to the ``third Kepler law'' for the angle $\ell_i$
\beq{Kep3}\ell_i\to \ell_i+\frac{{\eufm m}_i^3{\eufm M}_i^2}{\L_i^3}t\ ,\eeq
while all the remaining coordinates remain constant.

\nl
Let us recall that the change \equ{lin Del} can be extended to a four--dimensional canonical change 
\beq{2bintegration}(\L_i, \G_i, \ell_i, {\rm g}_i)\to ({\rm R}_i, \Phi_i, {\rm r}_i, \varphi_i)\qquad {\rm with}\quad \G_i=\Phi_i\eeq
involving the four coordinates $(\L_i, \G_i, \ell_i, {\rm g}_i)\in \real^2\times \torus^2$. 
Taking in count also the couples $({\rm H}_i, {\rm h}_i)$, one  lifts this map to a six--dimensional one
\beq{2bintegration**}(\L_i, \G_i, {\rm H}_i,  \ell_i, {\rm g}_i, {\rm h}_i)\to ({\rm R}_i, \Phi_i, {\rm H}_i, {\rm r}_i, \varphi_i, {\rm h}_i,)\qquad {\rm with}\quad \G_i=\Phi_i\eeq
via the identity on such couples.

\nl
In count of relations \equ{bar y bar x} and \equ{x***}, we correspondingly have two canonical maps
\beqa{Del*} \phi^\ppi_{\rm D}:&&\quad(\L_i, \G_i, {\rm H}_i, \ell_i, {\rm g}_i, {\rm h}_i)\to (y_{\rm D}^\ppi, x_{\rm D}^\ppi)\in \real^3\times \real^3\nonumber\\
\bar\phi_{\rm D}^\ppi:&&\quad (\L_i, \G_i, \ell_i, {\rm g}_i)\to (\bar y^\ppi_{\rm D}, \bar x^\ppi_{\rm D})\in \real^2\times \real^2 
\eeqa
which we shall refer  to as ``spatial'', ``planar'' Delaunay map, respectively. 
By similarity, we might call the map \equ{lin Del} ``linear'' Delaunay map.

\nl
The planar and the spatial maps  are related via
$$y^\ppi_{\rm D}(\L_i, \G_i, {\rm H}_i, \ell_i, {\rm g}_i, {\rm h}_i)={\rm R}^\ppt({\rm h}_i){\rm R}^\ppu(\iota_i)\bar y^\ppi_{\rm D}(\L_i, \G_i, \ell_i, {\rm g}_i)$$
$$x^\ppi_{\rm D}(\L_i, \G_i, {\rm H}_i, \ell_i, {\rm g}_i, {\rm h}_i)={\rm R}^\ppt({\rm h}_i){\rm R}^\ppu(\iota_i)\bar x^\ppi_{\rm D}(\L_i, \G_i, \ell_i, {\rm g}_i)$$
with $\cos\iota_i=\frac{{\rm H}_i}{\G_i}$.

\nl
We recall that relations $y^\ppi={\eufm m}_i\dot x^\ppi$  and \equ{Kep3} imply that $x^\ppi$ and $y^\ppi$ in \equ{Del*} (as well as $\bar x^\ppi_{\rm D}$ and $\bar y^\ppi_{\rm D}$, or ${\rm r}_i$ and ${\rm R}_i$) are related by
\beqa{yi}
y^\ppi(\L_i, \G_i, {\rm H}_i, \ell_i, {\rm g}_i, {\rm h}_i)&=&\frac{{\eufm m}_i^4{\eufm M}_i^2}{\L_i^3}\partial_{\ell_i} x^\ppi(\L_i, \G_i, {\rm H}_i, \ell_i, {\rm g}_i, {\rm h}_i)
\ .\eeqa

\nl
Let us now describe the meaning of the four variables at l. h. s. of \equ{2bintegration}.

\nl
As well known, the pull back orbits $ t\to (y^\ppi(t), x^\ppi(t))$  are ellipses ${\eufm E}_i$ with a focus in the origin. If $P^\ppi$, $a^\ppi$, $e^\ppi$ denote the perihelion, the semi--major axis, the eccentricity of ${\eufm E}_i$, $\cA^\ppi$ the area spanned from $P^\ppi$ to $x^\ppi$, $\cA^\ppi=\p (a^\ppi)^2\sqrt{1-(e^\ppi)^2}$ the total area of ${\eufm E}_i$, then $P^\ppi$ has the direction of ${\rm L}_i$ and the coordinates at left hand side in \equ{2bintegration} are defined by
$$
\arr{\dst
\L_i={\eufm m}_i\sqrt{{\eufm M}_i a^\ppi}\\ \\ \\
\G_i=\Phi_i=|{\rm C}^\ppi|
}\qquad \arr{\dst
\dst\ell_i=2\p\frac{\cA^\ppi}{\cA^\ppi_{\rm tot}}\\ \\
\dst{\rm g}_i=\a_{{\rm C}^\ppi}(n_i, P^\ppi)
}
$$
The angles $\ell_i$, ${\rm g}_i$ in this formula and ${\rm h}_i$ in \equ{reduct} are called ``mean anomaly'', ``argument of perihelion''and  ``longitude of the node''.

\nl
Delaunay coordinates are not defined when ${\rm C}^\ppi\parallel k^\ppt$ or the eccentricity of the $i^{\rm th}$ ellipse is zero. Then, following Poincar\'e,  one modifies them by considering 
the canonical variables
$$\L_i, \ul_i, \uh_i, \ux_i\ ,\up_i,\uq_i$$
where $\L_i$ are as above, and
$$\ul_i=\ell_i+{\rm g}_i+{\rm h}_i\quad \arr{\dst \uh_i=\sqrt{2(\L_i-\G_i)}\cos({\rm g}_i+{\rm h}_i)\\
\dst \ux_i=-\sqrt{2(\L_i-\G_i)}\sin({\rm g}_i+{\rm h}_i)}
\quad
\arr{\dst \up_i=\sqrt{2(\G_i-{\rm H}_i)}\cos({\rm h}_i)\\
\dst \uq_i=-\sqrt{2(\G_i-{\rm H})}\sin({\rm h}_i)\ .}
$$
One lets $\uz_i:=(\uh_i, \ux_i, \up_i, \uq_i)$ and $\uz=(\uz_1,\cdots, \uz_n)$.
The Poincar\'e maps
$$\phi^\ppi_{\rm P}:\qquad (\L_i, \ul_i, \uh_i, \ux_i\ ,\up_i,\uq_i)\to (y^\ppi, x^\ppi)$$
turn out to be regular also for $(\uh_i, \ux_i)=0$ or $(\up_i, \uq_i)=0$, which correspond to the vanishing of eccentricities and inclinations.
They have the form
\beq{poinc spatial}y^\ppi={\eufm R}^\ppi_{\rm P}(\L_i,\uz_i)\bar x^\ppi_{\rm P}(\L_i, \ul_i, \uh_i,\ux_i)\qquad x^\ppi={\eufm R}^\ppi_{\rm P}(\L_i,\uz_i)\bar x^\ppi_{\rm P}(\L_i, \ul_i, \uh_i,\ux_i)\eeq
where ${\eufm R}^\ppi_{\rm P}(\L_i,\uz_i)$ are unitary $3\times 3$ matrices depending only on $(\L_i,\uz_i)$ and reducing to the identity for $(\up_i, \uq_i)=0$ and 
\beq{poinc planar}\bar \phi^\ppi_{\rm P}:\qquad(\L_i, \ul_i, \uh_i,\ux_i)\to (\bar x^\ppi_{\rm P}(\L_i, \ul_i, \uh_i,\ux_i), \bar y^\ppi_{\rm P}(\L_i, \ul_i, \uh_i,\ux_i))\eeq
are the so--called planar Poincar\'e maps.
\subsection{Kepler maps}\label{Kepler maps}

\nl
Delaunay and Poincar\'e maps put the unperturbed term in \equ{helio} in the form
$$h_{\rm K}(\L)=-\sum_{1\le i\le n}\frac{{\eufm m}_i^3{\eufm M}_i^2}{2\L_i^2}\ .$$
Despite of this  good aspect, however, they  do not take in account the reduction of the integral ${\rm C}$ in \equ{C}. Namely, the transformed Hamiltonians\footnote{
If $\L=(\L_1, \cdots, \L_n)$, $\G=(\G_1, \cdots, \G_n)$, ${\rm H}=({\rm H}_1, \cdots,{\rm H}_n)$, $\ell=(\ell_1, \cdots, \ell_n)$, ${\rm g}=({\rm g}_1, \cdots, {\rm g}_n)$, ${\rm h}=({\rm h}_1, \cdots,{\rm h}_n)$ denote the collection of all Delaunay coordinates, we shall denote as
\beqano \phi_{\rm D}:\quad(\L, \G, {\rm H}, \ell, {\rm g}, {\rm h})\to (y, x)\in \real^3\times \real^3\qquad 1\le i\le n
\eeqano
the collection of all such maps. $\phi_{\rm P}$ is defined analogously.
}
$${\rm H}_{\rm D}:={\rm H}_{\rm hel}\circ\phi_{\rm D}:={\rm h}_{\rm K}(\L)+\m f_{\rm D}(\L,\G,{\rm H}, \ell, {\rm g}, {\rm h})$$
and 
$${\rm H}_{\rm P}:={\rm H}_{\rm hel}\circ\phi_{\rm P}:={\rm h}_{\rm K}(\L)+\m f_{\rm P}(\L,\ul,\uz)$$
have $3n$ degrees of freedom as well as ${\rm H}_{\rm hel}$. 

\nl
At this respect, we introduced in  \cite{pinzari14} the notion of ``Kepler map''. This is a generalization of Delaunay and Poincar\'e maps $\phi_{\rm D}$  and $\phi_{\rm P}$,  motivated by the fact that the proper degeneracies of $h_{\rm K}$ provide, as mentioned in the introduction, more freedom in the choice of canonical coordinates for the planetary system. 

\begin{itemize}
\item[\tiny\textbullet] Given $2n$   positive ``mass parameters'' ${\eufm m}_1$, $\cdots$, ${\eufm m}_n$, ${\eufm M}_1$, $\cdots$, ${\eufm M}_n$, a set ${\eufm X}\subset\real^{5n}$  and a bijection
\beqano
\t:\quad {\eufm X}&\to& \big\{({\eufm E}_1,\cdots,{\eufm E}_n)\in (E^3)^n\ ,\ {\eufm E}_i:\ {\rm ellipse}\big\}\nonumber\\
{\rm X}\in {\eufm X}&\to& \big({\eufm E}_1({\rm X}),\cdots,{\eufm E}_n({\rm X})\big)
\eeqano
which assigns to any ${\rm X}\in{\eufm X}$ an ordered set  of $n$ ellipses $({\eufm E}_1,\cdots,{\eufm E}_n)$  in the Euclidean space $E^3$ with strictly positive eccentricities and having a common focus ${\rm S}$, 
 we shall say that an {injective} map
\beqano
\phi:\quad ({\rm X}, \ell)\in {\cal D}^{6n}:={\eufm X}\times \torus^n\to (y_\phi({\rm X}, \ell),x_\phi({\rm X}, \ell))\in 
(\real^{3})^n\times (\real^{3})^n
\eeqano
is a {\sl Kepler map} if $\phi$
associates to  $({\rm X}, \ell)\in {\eufm X}\times\torus^n$, with $\ell=(\ell_1,\cdots,\ell_n)$ 
an element
$$(y_\phi({\rm X}, \ell),x_\phi({\rm X}, \ell))=(y_\phi^\ppu({\rm X}, \ell_1), \cdots, y_\phi^\ppn({\rm X},\ell_n), x_\phi^\ppu({\rm X}, \ell_1),\cdots, x_\phi^\ppn({\rm X}, \ell_n))$$
in the following way.
Letting, respectively, $P_\phi^\ppi({\rm X})$, $a_\phi^\ppi({\rm X})$, $e_\phi^\ppi({\rm X})$ and  $N_\phi^\ppi({\rm X})$   the   direction from ${\rm S}$ to the perihelion, the semi--major axis, the eccentricity and a prefixed direction of the plane of ${\eufm E}_i({\rm X})$, $x_\phi^\ppi({\rm X}, \ell_i)$ are the coordinates with respect to a prefixed orthonormal frame $(i,j,k)$ centered in ${\rm S}$ of the point of ${\eufm E}_i({\rm X})$ such that 
\beq{area}
\cA(x_\phi^\ppi({\rm X}, \ell_i))=\frac{1}{2}a^\ppi_\phi\sqrt{1-(e^\ppi_\phi)^2}\ell_i\quad \mod \p a^\ppi_\phi\sqrt{1-(e^\ppi_\phi)^2}
\eeq
 is the area spanned from $P_\phi^\ppi({\rm X})$ to $x_\phi^\ppi({\rm X}, \ell_i)$ relatively to the positive (counterclockwise) orientation determined by   $N_\phi^\ppi({\rm X})$ and 
\beq{yi***}
y_\phi^\ppi({\rm X}, \ell_i)={\eufm m}_i \sqrt{\frac{{\eufm M}_i}{(a^\ppi)^3}}\partial_{\ell_i} x^\ppi_\phi({\rm X}, \ell_i)\ .
\eeq

\item[\tiny\textbullet] One can consider {\sl canonical} Kepler maps, \ie, such that ${\rm X}\in{\eufm X}$ has the form ${\rm X}=({\rm P}, {\rm Q}, \L)$ where $\L=(\L_1,\cdots,\L_n)=({\eufm m}_1\sqrt{{\eufm M}_1a^\ppu_\phi}, \cdots, {\eufm m}_n\sqrt{{\eufm M}_na^\ppn_\phi})$, ${\rm P}=({\rm P}_1,\cdots,{\rm P}_{2n})$, ${\rm Q}=({\rm Q}_1,\cdots,{\rm Q}_{2n})$ and the map
$$(\L,\ell,{\rm P}, {\rm Q})\to (y,x)=(y^\ppu,\cdots,y^\ppn, x^\ppu, \cdots, x^\ppn)$$ preserves the standard 2-form:
$$
\sum_{1\le i\le n}d\L_i\wedge d\ell_i+\sum_{1\le i\le 2n}d{\rm P}_i\wedge d{\rm Q}_i=\sum_{1\le i\le n}\sum_{1\le j\le 3} dy^\ppi_j\wedge dx^\ppi_j$$
or, equivalently,
\beq{canonical map}\sum_{1\le i\le 2n}d{\rm P}_i\wedge d{\rm Q}_i=\sum_{1\le i\le n}d\G_i\wedge d{\rm g}_i+\sum_{1\le i\le n}d{\rm H}_i\wedge d{\rm h}_i\ .\eeq
\item[\tiny\textbullet] One can also consider canonical Kepler maps on manifolds, \ie, maps
$$(\L, \ell, \bar{\rm P}, \bar{\rm Q})\to (y, x)$$
where $\bar{\rm P}=(\bar{\rm P}_1,\cdots, \bar{\rm P}_m)$, $\bar{\rm Q}=(\bar{\rm Q}_1,\cdots, \bar{\rm Q}_m)$, with $m\le 2n$ and such that \equ{canonical map} holds with $\sum_{1\le i\le 2n}d{\rm P}_i\wedge d{\rm Q}_i$ replaced by $\sum_{1\le i\le m}d\bar{\rm P}_i\wedge d\bar{\rm Q}_i$.
\end{itemize}

\nl
Common properties to Kepler maps, with special attention to their application to the planetary hamiltonian (\eg, in the heliocentric form \equ{helio}) have been shortly discussed in \cite{pinzari14}. Here we just recall two facts.
\begin{itemize}
\item[\tiny\textbullet]  Relations \equ{area} and \equ{yi***} imply that the two--body Hamiltonian \equ{two bodies} becomes
$$h_{\rm K}^\ppi=-\frac{{\eufm m}_i{\eufm M}_i}{2a_\phi^\ppi}\ .$$
For a canonical map, one has just \equ{Kepler}.
\item[\tiny\textbullet] Due to \equ{yi***} and to 
equation
$$\sqrt{\frac{{\eufm M}_i}{(a^\ppi_\phi)^3}}\partial_{\ell_i}y_\phi^\ppi({\rm X}, \ell_i)=-{\eufm m}_i{\eufm M}_i\frac{x_\phi^\ppi({\rm X}, \ell_i)}{|x_\phi^\ppi({\rm X}, \ell_i)|^3}$$
(which follows from \equ{yi***} and the definition of $x_\phi^\ppi({\rm X}, \ell_i)$), one has 
\beq{zero average}\int_\torus \frac{x_\phi^\ppi({\rm X}, \ell_i)}{|x_\phi^\ppi({\rm X}, \ell_i)|^3} d\ell_i=0\qquad \int_\torus y_\phi^\ppi({\rm X}, \ell_i)d\ell_i=0\ .\eeq
Such relations give us the opportunity of noticing a dynamical equivalence between the 
heliocentric, Jacobi reduction, announced in \S \ref{Linear Momentum Reductions}, which is explained  in the following section.
\end{itemize}

\subsection{Dynamical equivalence of the heliocentric and Jacobi reduction}\label{Kepler maps*}

\nl
We just consider the case $n=2$. Let $\phi$ be a given Kepler map  in correspondence of mass parameters $({\eufm m}_i, {\eufm M}_i)_{1\le i\le 2}$
and $\tilde\phi$ the same Kepler map  in correspondence of  $(\tilde{\eufm m}_i, \tilde{\eufm M}_i)_{1\le i\le 2}$. Here, $({\eufm m}_i, {\eufm M}_i)_{1\le i\le 2}$ and $(\tilde{\eufm m}_i, \tilde{\eufm M}_i)_{1\le i\le 2}$ are as in \S \ref{Linear Momentum Reductions}.
Let  $f_{\rm hel, \phi}:=f_{\rm hel}\circ\phi$, $f_{\rm Jac, \tilde\phi}:=f_{\rm Jac}\circ\tilde\phi$ and put $$f_{\rm hel, \phi}^{\rm av}:=\frac{1}{(2\p)^2}\int_{\torus^2} f_{\rm hel, \phi}\,d\ell_1d\ell_2\qquad f_{\rm Jac, \tilde\phi}^{\rm av}:=\frac{1}{(2\p)^2}\int_{\torus^2} f_{\rm Jac, \tilde\phi}\,d\ell_1d\ell_2\ .$$ 

\nl
Then relations \equ{zero average} immediately imply\footnote{Recall \equ{Jac**} and ${\eufm m}_i=m_i+{\rm O}(\m)=\tilde{\eufm m}_i$, ${\eufm M}_i=m_0+{\rm O}(\m)=\tilde{\eufm M}_i$. Note that \equ{equivalence} holds also taking the integral only with respect to $\ell_2$.}
\beq{equivalence}f_{\rm hel, \phi}^{\rm av}= f_{\rm Jac, \tilde\phi}^{\rm av}+{\rm O}(\m)\ ,\eeq
with the common value of the two integrals being just the average of the Newtonian term, $\dst -\frac{m_1m_2}{(2\p)^2}\int_{\torus^2}\frac{d\ell_1d\ell_2}{|x^\ppu_\phi-x^\ppd_\phi|}+{\rm O}(\m)$. This relation is important in view of applications to KAM theory; for example, in relation to the proof of Arnold's Theorem, it says that the Birkhoff invariants that one finds using any of the two reductions differ\footnote{Observe the following remarkable consequence of this, noticed in \cite{abdullahA01}: while the first order Birkhoff invariants $\O_i$ associated to $f_{\rm hel, \phi}^{\rm av}$ satisfy {\sl identically} Herman resonance $\sum_i\O_i\equiv 0$, this resonance is instead verified only at the lowest order in $\m$ by the invariants associated to $f_{\rm Jac, \tilde\phi}^{\rm av}$.} just by ${\rm O}(\m)$.

\nl
Equality \equ{equivalence} has a dynamical explanation.

\nl
Averaging theory states that, for a properly--degenerate dynamical system
$${\rm H}(I,\varphi, u, v)=h(I)+\m f(I,\varphi, u, v)$$
there exists an  associated ``secular'' system
$${\rm H}^{\rm sec}(I,\varphi, u, v)=h(I)+\m f^{\rm av}(I, u, v)+{\rm O}(\m^2)$$
with $f^{\rm av}(I, u, v)=\frac{1}{(2\p)^n}\int_{\torus^n}f(I,\varphi, u, v)d\varphi$, related to ${\rm H}(I,\varphi, u, v)$ via a symplectic, close to the identity transformation. We recall that the transformation realizing this conjugation needs not to be unique, but  ${\rm H}^{\rm sec}$ is uniquely determined up to ${\rm O}(\m^2)$ if such transformation is chosen in the class of symplectic, $\m$--close to the identity ones.

\nl
Then one can conjugate the two Hamiltinians
$${\rm H}_{\rm hel, \phi}:={\rm H}_{\rm hel}\circ\phi=-\sum_{1\le i\le 2}\frac{{\eufm m}_i^3{\eufm M}_i^2}{2\L_i^2}+\m f_{\rm hel, \phi}(\L, \ell, P,Q)$$
and
$$\tilde{\rm H}_{\rm Jac, \tilde\phi}:={\rm H}_{\rm Jac}\circ\tilde\phi=-\sum_{1\le i\le 2}\frac{\tilde{\eufm m}_i^3\tilde{\eufm M}_i^2}{2\L_i^2}+\m f_{\rm Jac, \tilde\phi}(\L, \ell, P, Q)$$
 respectively to 
$${\rm H}_{\rm hel, \phi}^{\rm sec}=-\sum_{1\le i\le 2}\frac{{\eufm m}_i^3{\eufm M}_i^2}{2\L_i^2}+\m f^{\rm av}_{\rm hel, \phi}(\L,  P,Q)+{\rm O}(\m^2)$$
$$\tilde{\rm H}_{\rm Jac, \tilde\phi}^{\rm sec}=-\sum_{1\le i\le 2}\frac{\tilde{\eufm m}_i^3\tilde{\eufm M}_i^2}{2\L_i^2}+\m f_{\rm Jac, \tilde\phi}^{\rm av}(\L, P, Q)+{\rm O}(\m^2)\ .$$
Furthermore, relation \equ{hel/Jac}  
implies that ${\rm H}_{\rm hel, \phi}^{\rm sec}$ and $\tilde{\rm H}_{\rm Jac, \tilde\phi}^{\rm sec}$ are related by a symplectic, close to the identity transformation. Using finally\footnote{This follows from $${\eufm m}_1=\tilde{\eufm m}_1\ ,\quad {\eufm m}_2=\tilde{\eufm m}_2+{\rm O}(\m^2)\ ,\quad {\eufm m}_i{\eufm M}_i=\tilde{\eufm m}_i\tilde{\eufm M}_i=m_0m_i\ .$$
(The relation between ${\eufm m}_2$ and $\tilde{\eufm m}_2$ deserves to be remarked.)
}
$$
{\eufm m}_i^3{\eufm M}_i^2=\tilde{\eufm m}_i^3\tilde{\eufm M}_i^2+{\rm O}(\m^2)$$
and the uniqueness of the secular system associated to any of the two ${\rm H}_{\rm hel, \phi}$ and $\tilde{\rm H}_{\rm Jac, \tilde\phi}$ we have \equ{equivalence}.

\section{Examples of canonical Kepler maps}

\nl
The first classical example of Kepler map is clearly given by the Delaunay, Poincar\'e maps. In these cases, relations \equ{area} and \equ{yi***} are a consequence of the fact that such maps are canonical modifications of the six--dimensional map \equ{2bintegration**}.

\nl
On the other hand, in order that  \equ{area} and \equ{yi***}  are satisfied, it is not necessary to modify  \equ{2bintegration**}. Starting from the planar map \equ{2bintegration} or even  the linear one \equ{lin Del}   is sufficient.

\nl
In \S \ref{Jacobi Radau}--\S \ref{regularizations} we present examples of Kepler maps which are based on the planar map \equ{2bintegration}; in \S \ref{perihelia reduction} we provide an example based on the linear map  \equ{lin Del}.

\subsection{Jacobi--Radau's reduction of the nodes for three bodies}\label{Jacobi Radau}

\nl
The  ``Jacobi reduction of the nodes''  is a classical tool for reducing the number of degrees of freedom of \equ{helio} from six to four, in the case $n=2$. Its main properties are
\begin{itemize}
\item[\tiny\textbullet] it works only for the spatial problem (its planar limit is singular); 
\item[\tiny\textbullet]  uses the planar Delaunay maps \equ{2bintegration}; 
\item[\tiny\textbullet] 
 it can be regarded as an example of Kepler map on a eight--dimensional manifold (compare \equ{four manifold} below).
 \end{itemize}
 
 \nl
  Let us recall some history.

\nl
In 1842, Jacobi \cite{jacobi1842} discovered that  (after the linear momentum reduction) the twelve differential equations of the three--body problem dynamical system might  be reduced to a system of seven of equations: six of them of the first order, one of the second order. Essentially, four degrees of freedom. Clearly, what boils down is the reduction of the angular momentum integral. 
 
 \nl
Jacobi's procedure inherited a Hamiltonian aspect after the paper by Radau \cite{radau1868}, that now we explain. Let 
\beqa{3b}{\rm H}_{\rm 3b}({\rm R}, \Phi, {\rm H}, {\rm r}, \varphi, {\rm h})&=&\frac{{\rm R}_1^2}{2{\eufm m}_1}+\frac{\Phi_1^2}{2{\eufm m}_1{\rm r}_1^2}-\frac{{\eufm m}_1{\eufm M}_1}{{\rm r}_1}+\frac{{\rm R}_2^2}{2{\eufm m}_2}+\frac{\Phi_2^2}{2{\eufm m}_2{\rm r}_2^2}-\frac{{\eufm m}_2{\eufm M}_2}{{\rm r}_2}\nonumber\\
&+&\m \big(\frac{y^\ppu\cdot y^\ppd}{m_0}-\frac{m_1m_2}{|x^\ppu-x^\ppd|}\big)({\rm R}, \Phi, {\rm H}, {\rm r}, \varphi, {\rm h})\eeqa
be the heliocentric, three--body Hamiltonian  for $n=2$ planets, written in the coordinates \equ{reduct}. 

\nl
Recall that one has, for the heliocentric coordinates, the expressions in \equ{x***}. Such relations and   invariance by rotations of ${\rm H}_{\rm 3b}$ imply that ${\rm H}_{\rm 3b}$ depends upon the angles ${\rm h}_1$ and ${\rm h}_2$  only via the difference ${\rm h}_1-{\rm h}_2$. 
If we fix (as it is always possible to do) the initial frame ${\rm F}_0=(k^\ppu, k^\ppd, k^\ppt)$  such in a way that the total angular momentum ${\rm C}={\rm C}^\ppu+{\rm C}^\ppd$  coincides with the $k^\ppt$--axis, by definition, the nodes $n_1$, $n_2$ in \equ{D-nodes} are opposite one to the other, as it follows from
\beq{opposition of nodes}0=k^\ppt\times{\rm C}=k^\ppt\times{\rm C}^\ppu+k^\ppt\times{\rm C}^\ppd=n_1+n_2\ .\eeq
We then have 
\beq{nodes reduction}{\rm h}_2-{\rm h}_1=\p\qquad \mod 2\p\ .\eeq 
By this relation and again  invariance by rotations of ${\rm H}_{\rm 3b}$,   replacing Eq. \equ{x***} into ${\rm H}_{\rm 3b}$ \equ{3b} is equivalent to replace 
\beqa{xi}
&&x^\ppu={\rm R}_1(\iota_1) \bar x^\ppu({\rm R}_1, \Phi_1, {\rm r}_1, \varphi_1)\nonumber\\
&&x^\ppd={\rm R}_1(-\iota_2) \bar x^\ppd({\rm R}_2, \Phi_2, {\rm r}_2, \varphi_2)
\eeqa
 and, similarly, 
\beqa{yinew}
&&y^\ppu={\rm R}_1(\iota_1) \bar y^\ppu({\rm R}_1, \Phi_1, {\rm r}_1, \varphi_1)\nonumber\\
&&y^\ppd={\rm R}_1(-\iota_2) \bar y^\ppd({\rm R}_2, \Phi_2, {\rm r}_2, \varphi_2)\ .
\eeqa
Here we have used ${\rm R}^\ppt(\p){\rm R}^\ppu(\iota)={\rm R}^\ppu(-\iota){\rm R}^\ppt(\p)$ and an inessential shift of $\varphi_2$ by $\p$.

\nl
Fix
the  eight--dimensional\footnote{Equation $k^\ppt={\rm C}$ is equivalent to ${\rm C}_1={\rm C}_2=0$, hence corresponds to two conditions. Equation $k^\ppu=n_1$ also is bi-dimensional, since, in general, $k^\ppu$ varies among vectors verifying $k^\ppu\perp k^\ppt={\rm C}(\perp n_1)$. We then have four independent conditions in a 12-dimensional phase space.} manifold 
\beq{four manifold}\cM_{n_1}=\big\{k^\ppu=n_1\ ,\quad k^\ppt={\rm C}\big\}\ .\eeq  This corresponds to fix a rotating frame about ${\rm C}$
such that, with respect to it, node of the plane perpendicular to  ${\rm C}^\ppu$ is fixed. The aim is to find canonical coordinates for $\cM_{n_1}$. If we let\footnote{Note that, in principle, ${\rm G}$ is a function of $\Phi_1$, $\Phi_2$, ${\rm H}_1$, ${\rm H}_2$ and ${\rm h}_1-{\rm h}_2$.}
 ${\rm G}:=|{\rm C}|$, we may also write
the respective inclinations $\iota_1$ and $\iota_2$ as
$$\cos\iota_1=\frac{\Phi_1^2+{\rm G}^2-\Phi_2^2}{2\Phi_1{\rm G}}\qquad \cos\iota_2=\frac{\Phi_2^2+{\rm G}^2-\Phi_1^2}{2\Phi_2{\rm G}}\ .$$
These relations follow  considering the triangle formed by ${\rm C}^\ppu$, ${\rm C}^\ppd$ and ${\rm C}$. Similarly, the convex angle $\iota=\p-(\iota_1+\iota_2)$ formed by ${\rm C}^\ppu$ and ${\rm C}^\ppd$ is determined by
\beq{incli}\cos\iota=\frac{\Phi_1^2+\Phi_2^2-{\rm G}^2}{2\Phi_1\Phi_2}\ .\eeq
Again,  invariance by rotations of ${\rm H}_{\rm 3b}$ implies that only this angle is needed.

\nl
Let   \beqano
 {\rm H}_{\rm 3b, red}&=&\sum_{1\le i\le 2}\big(\frac{{\rm R}_i^2}{2{\eufm m}_i}+\frac{\Phi_i^2}{2{\eufm m}_i{\rm r}_i^2}-\frac{{\eufm m}_i{\eufm M}_i}{{\rm r}_i}\big)\nonumber\\
 &+&\m\big(\frac{y^\ppu\cdot y^\ppd}{m_0}-\frac{m_1m_2}{|x^\ppu-x^\ppd|}\big)({\rm R}, \Phi, {\rm r}, \phi,; {\rm G}) \eeqano
denote the
 four degrees of freedom function which is obtained from ${\rm H}_{\rm 3b}$ $({\rm R}$, $\Phi, {\rm H}$, ${\rm r}$, $\varphi, {\rm h})$ using relations \equ{bar y bar x}, \equ{xi}, \equ{yinew} and \equ{incli} with  ${\rm G}$ 
{\sl regarded  as an external fixed parameter}.  Following Arnold, we might say that the Hamiltonian
 ${\rm H}_{\rm 3b, red}({\rm R}, \Phi, {\rm r}, \varphi; {\rm G})$ ``resembles the Hamiltonian of  a certain planar problem\footnote{The verbatim citation is taken from \cite[end of \S 4, p. 141]{arnold63}. Here the reference to the planar problem is just relatively to the number of degrees of freedom. Note however that the spatial and the planar three--body problems  have the same number of degrees of freedom in the situation that the spatial problem is completely reduced and the planar one is not. The completely reduced planar problem has indeed three degrees of freedom, not four.}''. Radau proved (with a computational procedure) that the motion of the eight coordinates $({\rm R}_i, \Phi_i, {\rm r}_i, \varphi_i)_{1\le i\le 2}$ is governed by the Hamilton equations of ${\rm H}_{\rm 3b, red}({\rm R}, \Phi, {\rm r}, \varphi; {\rm G})$. This result  may be formulated by saying that, for any fixed value of ${\rm G}$, the imbedding
\beq{imb1}({\rm R}_1, {\rm R}_2, \Phi_1, \Phi_2, {\rm r}_1, {\rm r}_2, \varphi_1, \varphi_2)\in \real^{8}\to (y^\ppu, y^\ppd, x^\ppu, x^\ppd)\in \real^{12}\eeq
defined by Eqs. \equ{bar y bar x}, \equ{xi}, \equ{yinew} and \equ{incli}
is a canonical  map for the eight-dimensional manifold \equ{four manifold}. 

\nl
Let us incidentally mention that Radau's result (that he proved in about sixty pages) can be obtained as a consequence of a more general and global result that has been found one century later by A. Deprit \cite{deprit83}. Deprit also extended it to the case of more planets, as described in the next section.

\nl
To conclude this section dedicated to Jacobi--Radau's procedure, we describe two ``planetary'' modifications   of Jacobi--Radau variables, which are fitted to our needs.

\nl
Firstly,  we switch from the variables 
$({\rm R}_i, \Phi_i, {\rm r}_i, \varphi_i)$
to the variables on the left in \equ{2bintegration}. Then also the imbedding\footnote{The map \equ{imb2} is attributed to Jacobi by Arnold,  while the map \equ{imb1} is attributed to Radau by A. Deprit: compare \cite[beginning of \S 4, p.141]{arnold63} and \cite[end of p. 187]{deprit83}.}

\beq{imb2}(\L_1, \L_2, \G_1,\G_2, \ell_1, \ell_2, {\rm g}_1, {\rm g}_2)\in \real^{4}\times \torus^4\to (y^\ppu, y^\ppd, x^\ppu, x^\ppd)\in \real^{12}\eeq
is a canonical Kepler map, restricted  to the same manifold \equ{four manifold}. 
Such change reduces the two--body terms in \equ{3b} to their Keplerian form in \equ{Kepler}.

\nl
Secondly,   since  the map in \equ{2bintegration} is singular for  zero eccentricities,  it is customary to switch 
to a  set of Poincar\'e--like coordinates
\beq{Poinc after Jac}(\L_1, \L_2, \hat\ul_1,\hat\ul_2, \hat\uh_1,\hat\uh_2, \hat\ux_1,\hat\ux_2)
\eeq
where the ``secular variables'' $(\G_1,\G_2, {\rm g}_1, {\rm g}_2)$ have been replaced by 
$$\hat\uh_i=\sqrt{2(\L_i-\G_i)}\cos{\rm g}_i\qquad \hat\ux_i=-\sqrt{2(\L_i-\G_i)}\sin{\rm g}_i$$
 and the fast angles $(\ell_1, \ell_2)$ by
$$\hat\ul_i=\ell_i+{\rm g}_i\ .$$
The change
$$(\L_1, \L_2, \hat\ul_1,\hat\ul_2, \hat\uh_1,\hat\uh_2, \hat\ux_1,\hat\ux_2)\to (\L_1, \L_2, \G_1,\G_2, \ell_1, \ell_2, {\rm g}_1, {\rm g}_2)
$$
is canonical and  induces an imbedding
\beq{imb3}(\L_1, \L_2, \hat\ul_1,\hat\ul_2, \hat\uh_1,\hat\uh_2, \hat\ux_1,\hat\ux_2)\to 
(y^\ppu, y^\ppd, x^\ppu, x^\ppd)\in \real^{12}
\eeq
which
 regular\footnote{However, the imbedding \equ{imb3} is singular when the mutual inclination vanishes.  We should mention, at this respect, that, even though 
the change of coordinates \equ{imb3} is singular for zero mutual  inclination,  the Hamiltonian ${\rm H}_{\rm 3b, plt, reg}$ is instead regular. This nice aspect does not hold anymore when $n>2$.
} for $(\hat\uh_i, \hat\ux_i)=0$. Recall that $(\hat\uh_i, \hat\ux_i)=0$
  corresponds to the vanishing of the $i^{\rm th}$ eccentricity.

\nl
For future reference, we just mention that the maps in \equ{imb2} and \equ{imb3} have the form, respectively
\beqa{reduced 3 bodies action angle}
\arr{x^\ppu={\rm R}_1(\iota_1) \bar x_{\rm D}^\ppu(\L_1, \G_1,\ell_1,{\rm g}_1)\\
x^\ppd={\rm R}_1(-\iota_2) \bar x_{\rm D}^\ppd(\L_1, \G_1,\ell_1,{\rm g}_1)}\quad
\arr{y^\ppu={\rm R}_1(\iota_1) \bar y_{\rm D}^\ppu(\L_2, \G_2,\ell_2,{\rm g}_2)\\
y^\ppd={\rm R}_1(-\iota_2) \bar y_{\rm D}^\ppd(\L_2, \G_2,\ell_2,{\rm g}_2)
}
\eeqa
and
\beqa{reduced 3 bodies regularized}
\arr{x^\ppu={\rm R}_1(\iota_1) \bar x_{\rm P}^\ppu(\L_1, \hat\ul_1, \hat\uh_1,\hat\ux_1)\\
x^\ppd={\rm R}_1(-\iota_2) \bar x_{\rm P}^\ppd(\L_1, \hat\ul_1, \hat\uh_1,\hat\ux_1)}\quad
\arr{y^\ppu={\rm R}_1(\iota_1) \bar y_{\rm P}^\ppu(\L_2, \hat\ul_2, \hat\uh_2,\hat\ux_2)\\
y^\ppd={\rm R}_1(-\iota_2) \bar y_{\rm P}^\ppd(\L_2, \hat\ul_2, \hat\uh_2,\hat\ux_2)
}
\eeqa
where $(\L_i, \G_i, \ell_i,{\rm g}_i)\to(\bar x^\ppi_{\rm D}, \bar y^\ppi_{\rm D})$,  $(\L_i, \hat\ul_i, \hat\uh_i,\hat\ux_i)\to(\bar x^\ppi_{\rm P}, \bar y^\ppi_{\rm P})$ denote the planar  Delaunay, Poincar\'e maps, and $\iota_1$, $\iota_2$ have suitable expressions in such sets of coordinates.

\nl
Moreover,  we denote as
  \beqa{3b red}
 {\rm H}_{\rm 3b, red, plt}&=&-\sum_{1\le i\le 2}\frac{{\eufm m}_i^3{\eufm M}_i^2}{2\L_i}+\m\big(\frac{y^\ppu\cdot y^\ppd}{m_0}-\frac{m_1 m_2}{|x^\ppu-x^\ppd|}\big)(\L, \G, \ell, \g; {\rm G})\nonumber\\
 {\rm H}_{\rm 3b, red, plt,  reg}&=&-\sum_{1\le i\le 2}\frac{{\eufm m}_i^3{\eufm M}_i^2}{2\L_i}+\m\big(\frac{y^\ppu\cdot y^\ppd}{m_0}-\frac{m_1 m_2}{|x^\ppu-x^\ppd|}\big)(\L, \hat\ul, \hat\uz; {\rm G})
 \eeqa
 with $\hat\uz=(\hat\uh_1, \hat\uh_2,\hat\ux_1,\hat\ux_2)$,  the three--body written in the latter two different sets of variables.

\subsection{Boigey--Deprit's reduction for $(1+n)$ bodies}\label{Deprit' s reduction}
In 1893 A. Deprit \cite{deprit83}, strongly influenced by a previous paper by F. Boigey\footnote{The paper by Boigey deals with the case of four bodies ($n=3$). Boigey obtains a slightly less general result than the one by Deprit, since she constructs coordinates on a suitable manifold, not on the whole phase space.  Boigey's coordinates may be obtained as Deprit's coordinates for four bodies  restricted to Boigey's manifold.} \cite{boigey82} discovered a set of symplectic variables defined on the $6n$--dimensional space phase  that (suitably modified) ``unfold'' Jacobi--Radau's reduction of the nodes described in the previous section. 
Deprit's variables are defined in general for any many--particle system
$${\rm H}(y^\ppu, \cdots, y^\ppn, x^\ppu, \cdots, y^\ppn)$$
with general $n\ge 2$, which has the integral ${\rm C}$ in \equ{C} preserved. They are not specialized to the planetary problem, but may be adapted to it, as described in the next  \S \ref{regularizations}. 
Deprit's variables remained unnoticed for many years apart for some applications to the three--body problem  in which case, as said, they trivialize to Jacobi--Radau's reduction, described in \S \ref{Jacobi Radau}. See, for example, the paper \cite{ferrerO94}, where the Hamiltonian is just the same as in \cite{harrington69}, \cite{lidovZ76}, which are previous to Deprit's paper. Deprit's reduction has been next rediscovered (in the form given by the map \equ{Gab map}; see also \cite{chierchiaPi11a})  by the author during her PhD, under the motivation of their application to the direct proof of Arnold's  Theorem.

\nl
Analogously to the case of two planets, also Deprit's construction requires that certain inclinations appearing in the construction do not vanish. Let us start with the case $n=2$, in order to compare with  Jacobi--Radau's procedure.

\nl
Let ${\rm R}_1$, ${\rm R}_2$, ${\rm r_1}$, ${\rm r}_2$, $\Phi_1$ and $\Phi_2$ be as in \equ{reduct}. Let
$$\n_1:={\rm C}\times {\rm C}^\ppu\qquad \n_2:={\rm C}\times {\rm C}^\ppd\ .$$
Clearly, $\n_1$ and $\n_2$ are opposite, since
$$\n_1+\n_2={\rm C}\times {\rm C}^\ppu+{\rm C}\times {\rm C}^\ppd={\rm C}\times {\rm C}=0\ .$$
This relation is an extension of \equ{opposition of nodes}. Define 
\beq{n0}\n_0:=k^\ppt\times {\rm C}\ .\eeq
Then let
\beq{four variables}{\rm C}_3:={\rm C}\cdot k^\ppt\quad {\rm G}:=|{\rm C}|\quad \zeta:=\a_{k^\ppt}(k^\ppu, \n_0)\quad {\rm g}:=\a_{{\rm C}}(\n_0, \n_1)
\eeq
and
$$\phi_1:=\a_{{\rm C}^\ppu}(\n_1, x^\ppu)\quad \phi_2:=\a_{{\rm C}^\ppd}(\n_2, x^\ppd)=\a_{{\rm C}^\ppd}(\n_1, x^\ppd)+\p\ .$$

\nl
Deprit proved that the twelve--dimensional change
\beq{Dep map}({\rm C}_3, {\rm G}, {\rm R}_1, {\rm R_2}, \Phi_1, \Phi_2, \zeta, {\rm g}, {\rm r}_1, {\rm r}_2, \phi_1, \phi_2)\to (y^\ppu, y^\ppd, x^\ppu, x^\ppd)\eeq
is homogeneous. 

\nl
The main point (which holds also for the case $n\ge 3$ described below) is that ${\rm C}_3$, $\zeta$ and ${\rm g}$ are cyclic in any SO(3)--invariant Hamiltonian, since they are conjugated to integrals. The existence of two cyclic conjugated variables  is just the cause of the degeneracy of all the orders of  Birkhoff normal form   for the secular perturbing function written in Poincar\'e variables, found in \cite{chierchiaPi11c}.  The fact that an {\sl action}, ${\rm C}_3$, is cyclic is a remarkable situation, which generalizes what we already know from  the two-body Hamiltonian \equ{Kepler} (compare \S \ref{2 bodies}). We remark, at this respect that Deprit seemed not to notice this fact, since, at the end of \S 4, he underlines the cyclicality of $\zeta$ and ${\rm g}$ but not the one of ${\rm C}_3$.

\nl
Let us inspect the analytical expression of the map \equ{Dep map}. By the definitions, if $\iota_0$ is the convex angle formed by ${\rm C}$ and $k^\ppt$, determined by 
$$\cos\iota_0=\frac{\rm C_3}{{\rm G}}$$
then the map \equ{Dep map} 
is given by
\beqa{Dep map2b}
&&\arr{
\dst x^\ppu={\rm R}^\ppt(\zeta){\rm R}^\ppu(\iota_0){\rm R}^\ppt({\rm g}){\rm R}^\ppu(\iota_1)\bar x^\ppu\\ \\
\dst x^\ppd={\rm R}^\ppt(\zeta){\rm R}^\ppu(\iota_0){\rm R}^\ppt({\rm g}+\p){\rm R}^\ppu(\iota_2)\bar x^\ppd
}\nonumber\\\nonumber\\
&&\arr{
\dst y^\ppu={\rm R}^\ppt(\zeta){\rm R}^\ppu(\iota_0){\rm R}^\ppt({\rm g}){\rm R}^\ppu(\iota_1)\bar y^\ppu\\ \\
\dst y^\ppd={\rm R}^\ppt(\zeta){\rm R}^\ppu(\iota_0){\rm R}^\ppt({\rm g}+\p){\rm R}^\ppu(\iota_2)\bar y^\ppd
}\eeqa
with $\bar x^\ppu$, $\bar x^\ppd$, $\bar y^\ppu$, $\bar y^\ppd$ as in \equ{bar y bar x} with $\varphi_i$ replaced by $\phi_i$. The product ${\rm R}^\ppt(\zeta)$ ${\rm R}^\ppu(\iota_0)$ ${\rm R}^\ppt({\rm g})$ which is in common in all such formulae is negligible by  invariance by rotations. Therefore, in view of  expressing the three--body Hamiltonian \equ{3b} in terms of Deprit variables (at l. h. s. of \equ{Dep map}), the angles $\zeta$, $\iota_0$ and ${\rm g}$ can be fixed to any arbitrary value. Fixing them all to zero
reduces \equ{Dep map2b} to \equ{xi}--\equ{yinew} with $\varphi_i$ replaced by $\phi_i$. In particular, Radau's theorem mentioned in \S \ref{Jacobi Radau} is clarified\footnote{In \cite[end of p. 194]{deprit83} we read ``The intention [of this note] was to show how the global symmetry with respect to the group SO(3) ... affords a suitable coordinate system leading without artificiality the reduction of the nodes.''} from the geometrical point of view. Note in fact  that the manifold \equ{four manifold} in terms of the variables on the left in \equ{Dep map} has the simple expression
$${\rm C}_3={\rm G}_0\quad {\rm G}={\rm G}_0\quad \zeta=0\quad {\rm g}=0\ .$$

\nl
Let us now describe how Deprit  extended his reduction to the case of more than two planets.

\nl
Let ${\rm R}_i$, $\Phi_i$, ${\rm r}_i$ as in \equ{reduct}, with $1\le i\le n$.
 For $1\le j\le n$ define
\beq{Sj}
{\rm C}^\ppj:=x^\ppj\times y^\ppj\qquad {\rm S}^\ppj:=\sum_{j\le i\le n}{\rm C}^\ppi\eeq
so that ${\rm S}^\ppu={\rm C}$, ${\rm S}^\ppn={\rm C}^\ppn$. Then let $\n_0$ as in  \equ{n0} and, for $1\le j\le n-1$,
\beq{Deprit nodes}\n_j:={\rm S}^\ppj\times {\rm C}^{(j)}\quad 1\le j\le n-1\ ;\quad \n_n:={\rm S}^{(n-1)}\times {\rm C}^{(n)}
\ .\eeq
There are just $n$ independent $\n_j$'s (with $0\le j\le n-1$), since $\n_{n-1}$ and $\n_n$ are opposite. 
Assume that none of the $\n_j$'s vanishes. Then define
 ${\rm C}_3$, $\zeta$  as in \equ{four variables}  and\footnote{Here we use notations and definitions closer to Deprit's. In \cite{chierchiaPi11b} a different , but equivalent, definition of the variables $(\Psi, \psi)$ has been used. In \cite{deprit83} $({\rm C}_3, \zeta)$ are called $({\rm N}_0^*, \n_0^*)$, while $\Psi_0$, $\cdots$, $\Psi_{n-2}$, $\psi_0$, $\cdots$, $\psi_{n-2}$ are called $\Theta_0^*$, $\cdots$, $\Theta^*_{n-2}$, $\theta_0^*$, $\cdots$, $\theta^*_{n-2}$. Moreover $\Phi_1$, $\cdots$, $\Phi_n$, $\phi_1$, $\cdots$, $\phi_n$ are called $\Theta_1$, $\cdots$, $\Theta_n$, $\theta_1$, $\cdots$, $\theta_n$ in \cite{deprit83}.
 With respect to the notations used in \cite{chierchiaPi11b}, the correspondence is the following. Calling here $(\L^*$, $\G^*$, $\Psi^*$, $\ell^*$, $\g^*$, $\psi^*)$ the coordinates that in \cite{chierchiaPi11b} are named $(\L$, $\G$, $\Psi$, $\ell$, $\g$, $\psi)$, then $(\L_i$, $\G_i$, $\ell_i$, $\g_i)$=$(\L^*_{n-i+1}$, $\G^*_{n-i+1}$, $\ell^*_{n-i+1}$, $\g^*_{n-i+1})$,  and $(\Psi_j$, $\psi_j):=$ ($\Psi^*_{n-j-1}$, $\psi_{n-j-1}$, with $1\le i\le n$,  $-1\le j\le n-1$, $(\Psi_{-1}$, $\psi_{-1})$ as in \equ{new definitions}.  }
$$\arr{\dst \Psi_{j-1}:=|{\rm S}^\ppj|\\ \\
\dst \psi_{j-1}:=\a_{{\rm S}^\ppj}(\n_{j-1}, \n_j)
}\quad \quad \phi_i:=\a_{{\rm C}^\ppi}(\n_i, x^\ppi)$$
with $1\le i\le n$, $1\le j\le n-1$.
Then the variables $\Psi_0:={\rm G}$ and $\psi_0:={\rm g}$ are as in \equ{four variables}, with ${\rm C}$ the total angular momentum of $n$ planets. Define, finally, \beqano
&&{\rm R}:=({\rm R}_1,\cdots, {\rm R}_n)\qquad {\rm r}:=({\rm r}_1,\cdots, {\rm r}_n)\qquad \Phi:=(\Phi_1,\cdots, \Phi_n)\nonumber\\
&& \Psi:=({\rm C}_3, \Psi_0,\cdots, \Psi_{n-2})\quad \phi:=(\phi_1,\cdots, \phi_n)\quad \psi:=(\zeta, \psi_0, \cdots, \psi_{n-2})\ .
\eeqano
Deprit proved that the change
$$({\rm R}, \Phi, \Psi, {\rm r}, \phi, \psi)\in \real^{3n}\times \torus^{3n}\to (y^\ppu, \cdots, y^\ppn, x^\ppu, \cdots, x^\ppn)\in \real^{3n}\times \real^{3n}$$
is homogeneous.

\nl
The relative formulae to obtain any of the couples $(y^\ppi, x^\ppi)$ are a bit involved. They however can be found following  suitable\footnote{Here, ${\rm F}_0$ is the prefixed frame, ${\rm F}^*_1$ has the first axis in the direction of $\n_0$ and the third one in the direction of ${\rm S}^\ppu$; ${\rm F}_n:={\rm F}_n^*$, where ${\rm F}^*_j$, for $2\le j\le i\le n$ has the first axis in the direction of ${\rm S}^{(j-1)}\times {\rm S}^\ppj=-\n_{j-1}$ and the third one in the direction of ${\rm S}^\ppj$; ${\rm F}_i$, for $1\le i\le n-1$, has  the first axis in the direction of $\n_i={\rm S}^{(i-1)}\times {\rm C}^\ppj=\n_i$ and the third axis in the direction of ${\rm C}^\ppi$.} chains of frames
\beq{chain}{\rm F}_0\to{\rm F}_1^*\to\cdots\to {\rm F}_i^*\to {\rm F}_i\qquad 1\le i\le n\ ,\eeq
where ${\rm F}_0$ is a prefixed frame, ${\rm F}^*_1$, $\cdots$, ${\rm F}^*_i$ are auxiliary frames having their third axes in the direction of ${\rm S}^\ppu$, $\cdots$, ${\rm S}^\ppi$, respectively, ${\rm F}^\ppi$ is the $i^{\rm th}$ orbital frame, having its third axis in the direct of ${\rm C}^\ppi$.

\nl
Easy geometric arguments 
lead, for the map \equ{Dep map}, to expressions of the form
\beqa{Dep map*}
&&x^\ppi=\left\{\begin{array}{llll}
\dst \cR(\zeta, \iota_0)\cR(\psi_0,\iota_1)\bar x^\ppu\ &i=1\\ \\
\dst\cR(\zeta, \iota_0)\cR(\psi_0+\p,\iota_1^*)\cdots \cR(\psi_{i-2}+\p, \iota_{i}^*)\cR(\psi_{i-1}, \iota_i)\bar x^\ppi &2\le i\le n-1\\ \\
\dst\cR(\zeta, \iota_0)\cR(\psi_0+\p,\iota_1^*)\cdots \cR(\psi_{n-2}+\p, \iota_n^*)\bar x^\ppn &i=n
\end{array}
\right.\nonumber\\
&&y^\ppi=\left\{\begin{array}{llll}
\dst \cR(\zeta, \iota_0)\cR(\psi_0,\iota_1)\bar y^\ppu\ &i=1\\ \\
\dst\cR(\zeta, \iota_0)\cR(\psi_0+\p,\iota_1^*)\cdots \cR(\psi_{i-2}+\p, \iota_i^*)\cR(\psi_{i-1}, \iota_i)\bar y^\ppi &2\le i\le n-1\\ \\
\dst\cR(\zeta, \iota_0)\cR(\psi_0+\p,\iota_1^*)\cdots \cR(\psi_{n-2}+\p, \iota_n)\bar y^\ppn &i=n
\end{array}
\right.\nonumber\\
\eeqa
where $\bar x^\ppi$, $\bar y^\ppi$ depend on $({\rm R}_i, \Phi_i, {\rm r}_i, \phi_i)$ as in \equ{bar y bar x}, with $\varphi_i$ replaced by $\phi_i$; $\cR(\psi, \iota):={\rm R}^\ppt(\psi){\rm R}^\ppu(\iota)$; $\iota_i$ are the convex angles formed by ${\rm S}^\ppi$ and ${\rm C}^\ppi$, while $\iota^*_j$  are the convex angles formed by  ${\rm S}^\ppj$ and ${\rm S}^{(j+1)}$.
Such convex angles are determined by
$$\cos \iota_i=\frac{\Psi_{i-1}^2+\Phi_i^2-\Psi_i^2}{2\Psi_{i-1}\Phi_i}\qquad \cos \iota_j^*=\frac{\Psi_{j}^2+\Phi_j^2-\Psi_{j-1}^2}{2\Psi_{j}\Phi_j}$$
with $1\le i\le n$ and $1\le j\le n-1$. 

\nl
The description given in the formulae above (developed in \cite{pinzari-th09}; see also \cite[Appendix A]{chierchiaPi11b}) is different from Deprit's approach, who preferred to use quaternions and, at the end, did not provide complete formulae\footnote{In \cite[end of \S 4]{deprit83} we read ``{\sl The final expressions increase in complexity; it serves no purpose to enter the results in this Note.}''  }, since they seemed too complicate to him.
We should add, at this respect, that Deprit believed that his variables would never\footnote{At the end of p. 194 of Deprit's paper, we find ``{\sl Whether the new phase variables are practical in the General Theory of Perturbation is an open question. At least, for planetary theories, the answer is likely to be negative: the tree of kinetic frames imposes a recursive hierarchy without physical correspondence in the solar system.}''.} be useful.  His variables (in the rediscovered planetary version; compare the next section) have been instead the starting point for the proof of existence and non--degeneracy of the Birkhoff normal form for the planetary system and hence the complete, constructive proof of Arnold's Theorem \cite{pinzari-th09}, \cite{chierchiaPi11b}.

\nl
From the formulae \equ{Dep map*} we have an explicit explanation of the reduction: for SO(3)--invariant systems, the three first rotations by the angles $\zeta$, $\iota_0$ and $\psi_0$ in front of any $x^\ppi$ of $y^\ppi$ may be neglected and, as expected, only the variables
$$({\rm R}, \Phi, \hat\Psi, {\rm r}, \phi, \hat\psi)$$
where
$$\hat\Psi:=(\Psi_1,\cdots, \Psi_{n-2})\qquad \hat\psi:=(\psi_1,\cdots, \psi_{n-2})$$
 will appear and the integral/action $\Psi_0={\rm G}$  will play the r\^ole of an  external parameter.

\nl
For example, choosing the planetary heliocentric Hamiltonian ${\rm H}_{\rm hel}$ in \equ{helio}, one finds 
\beqa{Dep}{\rm H}_{\rm Dep}&=&\sum_{1\le i\le n}\big(\frac{{\rm R}_i^2}{2{\eufm m}_i}+\frac{\Phi_i^2}{2{\eufm m}_i{\rm r}_i^2}-\frac{{\eufm m}_i{\eufm M}_i}{{\rm r}_i}\big)\nonumber\\
&+&\m\sum_{1\le i<j\le n}\big(\frac{y^\ppi\cdot y^\ppj}{m_0}-\frac{m_i m_j}{|x^\ppi-x^\ppj|}\big)({\rm R}, \Phi, \hat\Psi, {\rm r}, \phi, \hat\psi; {\rm G})\ .\eeqa

\subsection{Planetary version of Deprit's variables and regularizations}\label{regularizations}
The planetary version of Deprit's coordinates (which is not mentioned in \cite{deprit83}) may obtained via the change  \beq{2bintegration*}(\L_i, \G_i, \ell_i, \g_i)\to ({\rm R}_i, \Phi_i, {\rm r}_i, \phi_i)\eeq defined via the integration of the two--body terms \equ{2bpc} terms appearing in  \equ{Dep}. With this change we then have a symplectic map
\beq{Gab map}(\L, \G, \Psi, \ell, \g, \psi)\to (y^\ppu, \cdots, y^\ppn, x^\ppu, \cdots, x^\ppn)\ ,\eeq
which turns out to be a Kepler map in the sense of \S \ref{Kepler maps}. More precisely, relations \equ{area}--\equ{yi***} are implied by the fact that the planar Delaunay map  is used.

\nl
The Hamiltonian  ${\rm H}_{\rm Dep}$ takes the ``planetary'' form
\beqno{\rm H}_{\rm Dep, plt}=-\sum_{1\le i\le n}\frac{{\eufm m}_i^3{\eufm M}_i^2}{2\L_i}+\m\sum_{1\le i<j\le n}\big(\frac{y^\ppi\cdot y^\ppj}{m_0}-\frac{m_i m_j}{|x^\ppi-x^\ppj|}\big)(\L, \G, \hat\Psi, \ell, \g,\hat\psi; {\rm G})\ .\eeqno
and coincides with the  Hamiltonian ${\rm H}_{\rm 3b, red, plt}$ in \equ{3b red} for $n=2$; for $n\ge 3$ it provides a suitable extension to it. It can be used for planetary theories in the case eccentricities and inclinations are not required to be small. Indeed, such occurrences are singular for ${\rm H}_{\rm Dep, plt}$: the singularities of inclinations is due to the definition of the variables $(\hat\Psi, \hat\psi)$ (compare \S \ref{Deprit' s reduction}); the singularities for eccentricities are introduced by the  maps \equ{2bintegration*}. For the reader who is  interested (for application purposes) to  the analytical expression of the map \equ{Gab map}, 
we just mention that such map is completely analogous to the map
 in \equ{Dep map*}, apart for taking  $\Phi_i=\G_i$ and replacing  $\bar x^\ppi({\rm R}_i, \Phi_i, {\rm r}_i, \phi_i)$, $\bar y^\ppi({\rm R}_i, \Phi_i, {\rm r}_i, \phi_i)$,  with $\bar x^\ppi_{\rm D}(\L_i, \G_i, \ell_i, \g_i)$, $\bar y^\ppi_{\rm D}(\L_i, \G_i, \ell_i, \g_i)$, where
 the maps $$\bar\phi^\ppi_{\rm D}:\qquad (\L_i, \G_i, \ell_i, \g_i)\to (\bar y^\ppi_{\rm D}(\L_i, \G_i, \ell_i, \g_i), \bar x^\ppi_{\rm D}(\L_i, \G_i, \ell_i, \g_i))$$
are as in \equ{Del*}. For $n=2$, it reduces to \equ{reduced 3 bodies action angle} neglecting the three first cyclic rotations by the angles $\zeta$, $\iota_0$, $\psi_0$.

\nl
To deal with the case when eccentricities or inclinations
may also vanish,  it is  possible to switch to new sets of coordinates, analogously to Poincar\'e's procedure for regularizing Delaunay coordinates. Unfortunately,  it is not possible to regularize {\sl all} singularities and, simultaneously, keep the number of degrees of freedom to $(3n-2)$.  This may be compared with the situation for the three--body case (see  \S \ref{Jacobi Radau}), where, as mentioned,  the coordinates \equ{Poinc after Jac} are suited for a four degrees of freedom  system,  but  the configuration with {\sl zero mutual inclination} is singular.

\nl
We may then choose if   regularizing {\sl all of} vanishing eccentricities or inclinations, at the cost of enhancing the number of degrees of freedom by one (``partial reduction''), or, alternatively,  {\sl all but one}  (``full reduction''). 
Moreover (see \S \ref{Full reduction and  symmetries})  full reduction, in general, breaks  down the natural symmetries of the problem; compare also \cite{maligeRL02}.

\subsubsection{Partial reduction (RPS variables)}\label{Partial reduction}
The complete regularization is obtained replacing the map \equ{Gab map} with a map
\beq{rps map}(\L,\l, z)\in \real^n\times \torus^n\times \real^{4n}\to (y^\ppu, \cdots, y^\ppn, x^\ppu, \cdots, x^\ppn)\eeq
where
$$z=(\eta_1,\cdots, \eta_n, \xi_1,\cdots, \xi_n, p_0,\cdots, p_{n-1}, q_0,\cdots, q_{n-1})$$
which is regular for $(\eta_i, \xi_i)$, $(p_j, q_j)$ approaching $(0,0)$ ($1\le i\le n$, $0\le j\le n-1$) and such that $(\eta_i, \xi_i)=0$ corresponds to the vanishing of the $i^{\rm th}$ eccentricity, $(p_j, q_j)=0$ corresponds to the vanishing of the node $\n_j$ defined in \S \ref{Deprit' s reduction}. The explicit  formulae of  the map
\equ{rps map}
 are given in \cite{pinzari-th09}, \cite{chierchiaPi11b}. The variables $(\L,\l,z)$ have been found in \cite{pinzari-th09} and have been named ``RPS'' (Regular, Planetary and Symplectic) in \cite{chierchiaPi11b}. 
 Their complete definition is
 \beqano&&\l_j=\ell_i+\g_i+\psi^{(i-1)}\ ,
\qquad{\rm where}\qquad \psi^\ppj=\sum_{-1\le k\le j} \psi_k\ ,
\nonumber  \\
&&\arr{\dst \eta_i=\sqrt{2(\L_i-\G_i)}\cos\big(\g_i+\psi^{(i-1)}\big)\\ \\
\dst \xi_i=-\sqrt{2(\L_i-\G_i)}\sin\big(\g_i+\psi^{(i-1)}\big)
}\nonumber\\
&&
\arr{\dst p_{j}=\sqrt{2(\G_{j}+\Psi_{j}-\Psi_{j-1})}\cos\psi^{(j-1)}\\ \\
\dst q_{j}=-\sqrt{2(\G_{j+1}+\Psi_{j+1}-\Psi_j)}\sin\psi^{(j-1)}
}
\eeqano
with
\beq{new definitions}\psi_{-1}:=\zeta\ , \ \Psi_{-1}:={\rm C}_3\ ,\ \G_0:=0\ ,\ 1\le i\le n\ ,\ 0\le j\le n-1\ .\eeq

 \nl
 Main points are
\begin{itemize}
\item[\tiny\textbullet]
The couple\footnote{Beware that here we are using different notations with respect to
 \cite{chierchiaPi11b}. Letting $\L^*$, $\l^*$ $\eta^*$, $\xi^*$, $p^*$, $q^*$ the variables named $\L$, $\l$ $\eta$, $\xi$, $p$, $q$ in  \cite{chierchiaPi11b}, the correspondence is: $(\L_i$, $\l_i$, $\eta_i$, $\xi_i)$=$(\L^*_{n-i+1}$, $\l^*_{n-i+1}$, $\eta^*_{n-i+1}$, $\xi^*_{n-i+1})$; $(p_j$, $q_j)=$ $(p^*_{n-j}$, $q^*_{n-j})$, with $1\le i\le n$, $0\le j\le n-1$. }  $(p_0, q_0)$ is cyclic and plays the r\^ole of $({\rm C}_3, \zeta)$ in the set on the left in \equ{Gab map}. Indeed,
$$p_0=\sqrt{2({\rm G}-{\rm C}_3)}\cos\zeta\qquad q_0=-\sqrt{2({\rm G}-{\rm C}_3)}\sin\zeta$$
are both integrals;
\item[\tiny\textbullet] The variables $(\L,\l, z)$ behave exactly as the Poincar\'e variables $(\L,\ul,\uz)$ of \S \ref{2 bodies}  for what concerns D'Alembert rules\footnote{D'Alembert rules are the expressions of rotation and reflection transformations in terms of  Poincar\'e variables. See, \eg,  \cite{chierchiaPi11b}, \cite{chierchiaPi11c}, \cite{pinzari14} for notices.}. In particular, $\bar z=0$, where $$\bar z:=(\eta, \xi, \bar p, \bar q):=(\eta_1,\cdots, \eta_n, \xi_1,\cdots, \xi_n, p_1,\cdots, p_{n-1}, q_1,\cdots, q_{n-1})$$
is $z$ deprived of $(p_0, q_0)$, is an elliptic equilibrium point for the secular system associated to ${\rm H}_{\rm rps}$, where ${\rm H}_{\rm rps}(\L,\l,\bar z)$ is ${\rm H}_{\rm hel}$ in \equ{helio} expressed in the variables $(\L,\l,z)$. More details are in \cite{pinzari-th09}, \cite{chierchiaPi11b}, \cite{chierchiaPi11c}.
\item[\tiny\textbullet] The formulae relating the change of coordinates \equ{rps map}
have the form 
$$x^\ppi=\left\{
\begin{array}{llll}
\dst {\eufm R}_0^*{\eufm R}_1 \bar x^\ppu_{\rm P}\quad &i=1\\ \\
\dst{\eufm R}_0^*\cdots{\eufm R}_i^*{\eufm R}_i \bar x^\ppi_{\rm P}&2\le i\le n
\end{array}
\right.
\quad y^\ppi=\left\{
\begin{array}{llll}
\dst {\eufm R}_0^*{\eufm R}_1 \bar y^\ppu_{\rm P}\quad &i=1\\ \\
\dst{\eufm R}_0^*\cdots{\eufm R}_i^*{\eufm R}_i \bar y^\ppi_{\rm P}&2\le i\le n
\end{array}
\right.$$
where ${\eufm R}_n^*=\id$, ${\eufm R}_0^*$, $\cdots$, ${\eufm R}_i^*$, ${\eufm R}_i$ are regular, unitary $3\times 3$ matrices depending only on $(\L,z)$; ${\eufm R}_0^*$ is completely negligible, since reduces to the identity for $(p_0, q_0)=0$;  the remaining matrices 
${\eufm R}_1^*$, $\cdots$, ${\eufm R}_i^*$, ${\eufm R}_i$ depend only on $(\L, \bar z)$ and
 reduce to the identity for $(\bar p, \bar q)=0$. In particular, as well as the spatial Poincar\'e maps \equ{poinc spatial}, also the map \equ{rps map}  reduces to the planar Poincar\'e maps \equ{poinc planar} for $(p, q)=0$.
 \item[\tiny\textbullet]  Even in the $n=2$ case, the  map \equ{rps map} is qualitatively different from the map \equ{reduced 3 bodies regularized} obtained via Jacobi reduction of the nodes. Indeed,  the map  \equ{rps map}  is regular for all vanishing eccentricities and inclinations  and  reduces the system to five d.o.f., while the map  \equ{reduced 3 bodies regularized} is regular for all vanishing eccentricities but singular for vanishing mutual inclination and reduces the system to   four d.o.f. The natural extension of the map \equ{reduced 3 bodies regularized} to the case $n\ge 3$ is described in the following section.
 \end{itemize}
\subsubsection{Full reduction}\label{Full reduction} 
Another full reduction (besides the one described in \S \ref{Deprit' s reduction}) is obtained giving up the regularization of the singular manifold $ \{\n_1=0\}$ (where $\n_1$ is defined in \equ{Deprit nodes}), corresponding of the parallelism of the three vectors ${\rm C}$, ${\rm C}^\ppu$ and $ {\rm S}^\ppd$. The remaining
eccentricities and  inclinations may be described instead with regularity.
The map realizing this reduction, restricted to the rotating manifold \equ{four manifold}, has the form of an imbedding
$$(\L,\hat\l, \hat z)\in \real^n\times \torus^n\times \real^{4n-4}\to (y^\ppu, \cdots, y^\ppn, x^\ppu, \cdots, x^\ppn)$$
with
\beqa{Rij}
&&x^\ppi=\left\{
\begin{array}{llll}
\dst
{\rm R}^\ppu(\iota_1)
\bar x^\ppu_{\rm P}\quad &i=1\\ \\
\dst
{\rm R}^\ppu(-\iota_1^*)
{\eufm R}_2\cdots{\eufm R}_i^*{\eufm R}_i \bar x^\ppi_{\rm P}&2\le i\le n
\end{array}
\right.
\nonumber\\
&&y^\ppi=\left\{
\begin{array}{llll}
\dst
{\rm R}^\ppu(\iota_1)
\bar y^\ppu_{\rm P}\quad &i=1\\ \\
\dst
{\rm R}^\ppu(-\iota_1^*)
{\eufm R}_2\cdots{\eufm R}_i^*{\eufm R}_i \bar y^\ppi_{\rm P}&2\le i\le n
\end{array}
\right.\eeqa
with $\psi_0^*:=\psi_0+\zeta$ and ${\eufm R}_i^*$, ${\eufm R}_j$, $\bar x^\ppj_{\rm P}$, $\bar y^\ppj_{\rm P}$  as in the previous section. The formulae \equ{Rij} naturally extend the formulae for the $n=2$ case in
 \equ{reduced 3 bodies regularized}.
 
 \nl
 We shall denote as
 \beqno{\rm H}_{\rm rps}=-\sum_{1\le i\le n}\frac{{\eufm m}_i^3{\eufm M}_i^2}{2\L_i}+\m\sum_{1\le i<j\le n}\big(\frac{y^\ppi\cdot y^\ppj}{m_0}-\frac{m_i m_j}{|x^\ppi-x^\ppj|}\big)(\L,\l, \bar z)\eeqno
 and
 \beqno{\rm H}_{\rm full\ red}=-\sum_{1\le i\le n}\frac{{\eufm m}_i^3{\eufm M}_i^2}{2\L_i}+\m\sum_{1\le i<j\le n}\big(\frac{y^\ppi\cdot y^\ppj}{m_0}-\frac{m_i m_j}{|x^\ppi-x^\ppj|}\big)(\L,\hat\l, \hat z;{\rm G})\eeqno
 the planetary system, written in the variables described in \S \ref{Partial reduction}--\S\ref{Full reduction} respectively.

 \subsection{Full reduction and reflection symmetries}\label{Full reduction and  symmetries}

Reducing completely all the integrals implies that symmetries related to them will be eliminated from the Hamiltonian. 

\nl
A natural question is whether it is possible to reduce the system by the integrals and keep, on the other side, parities due to those  symmetries  which are not related to integrals; like, for example, the transformations \equ{reflections}.

\nl
To begin with, let us investigate what happens to reflection symmetries when the system is completely reduced by rotations according to the reductions for $n=2$, $n\ge 3$, described, respectively, in   \S \ref{Jacobi Radau} or \S \ref{Full reduction}.

\nl
In the three--body case we have that the system retains a suitable parity in the  variables $\hat\uz=(\hat\uh, \hat\ux)$, due to reflections. This can be seen analyzing the the maps \equ{reduced 3 bodies regularized}:
the changes
\beqa{reflectionsP}
&&(\L_i, \hat\ul_i, \hat\uh_i, \hat\ux_i)\to (\L_i, -\hat\ul_i, -\hat\uh_i, \hat\ux_i)\nonumber\\
&&(\L_i, \hat\ul_i, \hat\uh_i, \hat\ux_i)\to (\L_i, \p-\hat\ul_i, \hat\uh_i, -\hat\ux_i)
\eeqa
correspond, respectively, to transform the projections $(\bar y^\ppi_{\rm P}, \bar x^\ppi_{\rm P})=\big((\bar y^\ppi_{\rm P, 1}, \bar y^\ppi_{\rm P, 2}), (\bar x^\ppu_{\rm P, 1}, \bar x^\ppi_{\rm P, 2})\big)$ of the planar Delaunay map  as follows:
$$(\bar y^\ppi_{\rm P, 1}, \bar y^\ppi_{\rm P, 2}), (\bar x^\ppi_{\rm P, 1}, \bar x^\ppi_{\rm P, 2})\to (-\bar y^\ppi_{\rm P, 1}, \bar y^\ppi_{\rm P, 2}), (\bar x^\ppi_{\rm P, 1}, -\bar x^\ppi_{\rm P, 2})$$
$$(\bar y^\ppi_{\rm P, 1}, \bar y^\ppi_{\rm P, 2}), (\bar x^\ppi_{\rm P, 1}, \bar x^\ppi_{\rm P, 2})\to (\bar y^\ppi_{\rm P, 1}, -\bar y^\ppi_{\rm P, 2}), (-\bar x^\ppi_{\rm P, 1}, \bar x^\ppi_{\rm P, 2})\ .$$
In view of the fact that the angles $\iota_1$, $\iota_2$ remain unchanged under \equ{reflectionsP}, we shall have that the coordinates  $(y^\ppi, x^\ppi)=\big((y^\ppi_1, y^\ppi_2, y^\ppi_3), (x^\ppi_1, x^\ppi_2, x^\ppi_3)\big)$ into the formulae \equ{reduced 3 bodies regularized}  undergo the following transformations
\beqa{yx}
&&( y^\ppi_{1},  y^\ppi_{2}, y^\ppi_{3}), ( x^\ppi_{1},  x^\ppi_{2}, x^\ppi_{3})\to (- y^\ppi_{1},  y^\ppi_{2}, y^\ppi_{3}), ( x^\ppi_{1}, - x^\ppi_{2}, -x^\ppi_{3})\nonumber\\
&&( y^\ppi_{1},  y^\ppi_{2}, y^\ppi_{3}), ( x^\ppi_{1},  x^\ppi_{2}, x^\ppi_{3})\to ( y^\ppi_{1}, - y^\ppi_{2}, -y^\ppi_{3}), (- x^\ppi_{1},  x^\ppi_{2}, x^\ppi_{3})
\eeqa
and hence ${\rm H}_{\rm 3b, red, reg}$ is left unvaried. We then have that the point $\hat\uz=0$ is  an equilibrium point for the averaged secular perturbation, still in the case the reduction is performed completely. This equilibrium turns out to be elliptic in the case of small mutual inclination; it is hyperbolic for large inclinations. In correspondence of such situations, unstable  tori with two frequencies for large inclinations have been found in \cite{jefferysM68}; stable  tori with four (maximal) frequencies have been found in \cite{robutel95},  when the inclinations are small.

\nl
The situation is rather different when $n\ge 3$. In that case, we have to face up formulae \equ{Rij}. Such formulae have a more complicate structure than \equ{reduced 3 bodies regularized}. It is possible to see that  \equ{yx}  holds for $i=n-1$, $n$, but not in general. And in fact it turns out that the Taylor expansion of secular perturbing function around $\hat z=0$ for the fully reduced system described in \S \ref{Full reduction}  contains also powers with odd degree. The point $\hat z=0$ is no longer an equilibrium and the construction of the Birkhoff normal form is prevented.

\nl
This problem has been treated {\sl locally} in \cite{chierchiaPi11b} interchanging the order of operations: postponing the full reduction after the Birkhoff--normalization of the partially reduced system and then applying an Implicit Function Theorem procedure definitely restores (in the range of small mutual inclinations) a suitable Birkhoff normal form (in particular, the elliptic equilibrium) also for the completely reduced system. The  cost is that  the reduction is defined only in the local domain of the Birkhoff transformation.
\subsection{``Perihelia reduction'': A symmetric full reduction of SO(3) invariance}\label{perihelia reduction}

In this section we present  Kepler map for the planetary problem which reduces the number of degrees of freedom to $(3n-2)$ and, simultaneously, keeps memory of reflection invariance for any $n\ge 2$.

\nl
This  map
\beq{P*}(\phi_{\rm P_*})^{-1}:\qquad (y^\ppu, \cdots, y^\ppn, x^\ppu, \cdots, x^\ppn)\to {\rm P}_*=(\L, \chi,\Theta,\ell,\k,\vartheta)\eeq
where
\beqano
\begin{array}{lll}
&\dst\L=(\L_1,\cdots,\L_n)\in \real^n\quad &\ell=(\ell_1,\cdots,\ell_n)\in \torus^n\\\\
&\dst\chi=(\chi_0,\bar\chi)\in \real\times \real^{n-1}& \k=(\k_0,\bar\k)\in \torus\times\torus^{n-1}\\\\
&\dst\Theta=(\Theta_0,\bar\Theta)\in \real\times \real^{n-1}&\vartheta=(\vartheta_0,\bar\vartheta)\in \torus\times\torus^{n-1}
\end{array}
\eeqano
with 
\beqano
&&\bar\chi=(\chi_1,\cdots,\chi_{n-1})\ ,\quad \bar\k=(\k_1,\cdots,\k_{n-1})\nonumber\\\nonumber\\
&& \bar\Theta=(\Theta_1,\cdots,\Theta_{n-1})\ ,\quad \bar\vartheta=(\vartheta_1,\cdots,\vartheta_{n-1})
\eeqano
is defined as follows.

\vskip.1in
\noi 
Let, as in \S \ref{Kepler maps},  ${\eufm m}_1$, $\cdots$, ${\eufm m}_n$, ${\eufm M}_1$, $\cdots$, ${\eufm M}_n$ be fixed mass parameters; 
let   $a^\ppi\in \real_+$,  $e^\ppi$, $P^\ppi\in \real^3$, with $|P^\ppi|=1$, denote  the semi--major axis, eccentricity, the   {direction of the perihelion} of the $i^{\rm th}$ instantaneous ellipse ${\eufm E}_i={\eufm E}_i(y^\ppi, x^\ppi)$ generated by the two--body Hamiltonian $h^\ppi_{\rm 2b}$ in \equ{two bodies} with initial datum $(y^\ppi, x^\ppi)$; let 
 $\cA^\ppi$ the  area spanned by $x^\ppi$ on ${\eufm E}_i$ with respect to $P^\ppi$. Let  and ${\rm C}^\ppi$, ${\rm S}^\ppi$ be as in \equ{Sj}.
Define, finally, the following $n$ couples of nodes, $(\widetilde\n_j,\widetilde{\rm n}_j)_{1\le j\le n}$

\beq{good nodes}\widetilde\n_1:=k^\ppt\times {\rm C}\ ,\quad \widetilde{\rm n}_j:={\rm S}^\ppj\times P^\ppj\ ,\quad \widetilde\n_{j+1}:=P^{(j)}\times {\rm S}^{(j+1)}\ ,\quad \widetilde{\rm n}_n:=P^\ppn\eeq
with $1\le j\le n-1$. Assume that such nodes do not vanish. Then let

\beqa{belle*}
\begin{array}{llllrrr}
\dst \Theta_{j-1}=\left\{
\begin{array}{lrrr}
\dst{\rm C}_3:={\rm C}\cdot k^\ppt
\\
\\
\dst
{\rm S}^{(j)}\cdot P^{(j-1)}
\end{array}
\right.& \vartheta_{j-1}=\left\{
\begin{array}{lrrr}
\dst\zeta:=\a_{k^\ppt}(k^\ppu, \widetilde\n_1)\qquad& j=1\\
\\
\dst \a_{P^{(j-1)}}(\widetilde{\rm n}_{j-1}, \widetilde\n_{j})&2\le j\le n
\end{array}
\right.\\ \\
\dst\chi_{j-1}:=\left\{
\begin{array}{lrrr}{\rm G}=|{\rm S}^\ppu|
\\
\\
|{\rm S}^\ppj|
\end{array}
\right.
&
\k_{j-1}:=\left\{
\begin{array}{lrrr}{\eufm g}:=\a_{{\rm S}^\ppu}(\widetilde\n_1, \widetilde{\rm n}_1)\qquad& j=1\\
\\
\a_{{\rm S}^\ppj}(\widetilde\n_j, \widetilde{\rm n}_j)&2\le j\le n
\end{array}
\right.
\\ \\
\L_i:={\eufm M}_i\sqrt{{\eufm m}_i a^\ppi}\qquad & \ell_i:=2\p\frac{\cA^\ppi}{\cA_{\rm tot}^\ppi}:= {\rm mean\ anomaly\ of}\ x^\ppi 
\ {\rm on}\ {\eufm E}_i
\end{array}
\eeqa
The map \equ{P*} is Kepler map in the sense of \S \ref{Kepler maps}. Relations \equ{area}--\equ{yi***}  hold by definition of the instantaneous ellipses ${\eufm E}_i(y^\ppi, x^\ppi)$.

\nl
Main points are
\begin{itemize}

\item[\tiny\textbullet]  The reduction realized by the variables \equ{P*} is based  on a new chain of frames 

\beq{chain2}{\rm F}_0\to \widetilde{\rm F}_1\to {\rm G}_1\to\cdots\to\widetilde{\rm F}_n\to {\rm G}_n\ .\eeq
This chain  contains $(2n)$ changes. This may be compared with the chains \equ{chain}, which instead may  contain up to $n$ changes. In the chain  \equ{chain2},  ${\rm F}_0$ is a prefixed initial frame;   $\widetilde{\rm F}_j$, with $1\le j\le n$ are frames which (analogously to the frames  ${\rm F}_j$ in \equ{chain}) have  their respective third axes  directed towards ${\rm S}^\ppj$.  The  frames  ${\rm G}_j$ are completely new (and this is the main difference with Deprit's reduction): they  have  their third axes are directed towards the perihelia $P^\ppj$ of instantaneous ellipses. For this reason we call this reduction as ``Perihelia reduction'', while Deprit's reduction might be named ``inclinations reduction''. Clearly, also the inclinations play a r\^ole  in \equ{chain2}, by means of the $\widetilde{\rm F}_j$'s, but the main point is that  the orbital frames, with their third axes in the direction of ${\rm C}^\ppj$ (corresponding to the ${\rm F}_j$ in \equ{chain}), are not part of \equ{chain2}. In other words, ${\rm C}^\ppj$  are not independent vectors in the reduction (while $P^\ppj$  are so) and their lengths $\G_j=|{\rm C}^\ppj|$ are not independent actions. 

\item[\tiny\textbullet] 
The   coordinates  $(\L,\bar\chi, \bar\Theta, \ell, \bar\k,\bar\vartheta)\in \real^n\times \real^{n-1}\times \real^{n-1}\times \torus^n\times \torus^{n-1}\times \torus^{n-1}$ are canonical coordinates  for $(6n-4)$--dimensional manifold
$$\cM_{\tilde\n_1}=\{k^\ppu=\tilde\n_1\qquad k^\ppt={\rm C}\}\ .$$
Such manifold corresponds to fix a rotating frame with the third axis parallel to ${\rm C}$ and such that, with respect to it, the node $\tilde\n_1$ determined by the plane orthogonal to the first perihelion $P^\ppu$ is fixed. This  should be compared with the corresponding manifold \equ{four manifold} of Deprit's reduction.

\item[\tiny\textbullet]  Singularities of the coordinates \equ{belle*} appear in correspondence of vanishing eccentricities (in which case the perihelia $P^\ppi$ are not defined) or cases of parallelism of two consecutive frames in the chain \equ{chain2}. The parallelism between two frames $\widetilde{\rm F}_j$, $\widetilde{\rm F}_{j+1}$ is not a singularity. In particular, the coordinates \equ{belle*} are well defined also in the case of planar limit (which corresponds to all the $\widetilde{\rm F}_j$'s parallel one with another); contrarily to what happens for Jacobi--Radau and Boigey--Deprit's reductions.

\item[\tiny\textbullet] The fact that the perihelia $P^\ppi$ are independent directions allows for a symmetry by reflections. (In the reductions by Boigey--Deprit, the directions of the perihelia are not independent, since they are constrained to be  orthogonal to the ${\rm C}^\ppj$'s, which are independent directions.) Transformations
$$(\L,\chi, \Theta, \ell, \k, \vartheta)\to (\L,\chi, -\Theta. \ell, \k, 2k\p-\vartheta)\qquad k\in \integer^n$$
correspond to 
$$\big((y^\ppi_1, y^\ppi_2, y^\ppi_3), (x^\ppi_1, x^\ppi_2, x^\ppi_3)\big)\to \big((y^\ppi_1, -y^\ppi_2, y^\ppi_3), (x^\ppi_1, -x^\ppi_2, x^\ppi_3)\big)\ .$$

\end{itemize}
The coordinates \equ{belle*} have been  presented in the note \cite{pinzari14}, together with a sketchy illustration of an application of them. The proof of their canonical character (and of the result announced in that note) is deferred to a subsequent paper.


%
%
\addcontentsline{toc}{section}{References}
\def\cprime{$'$}

\end{document}